\documentclass{amsart}
\usepackage{bbm,amsmath,amsfonts,amssymb}
\usepackage{tikz}
\RequirePackage[colorlinks,citecolor=blue,urlcolor=blue]{hyperref}
\newtheorem{theorem}{Theorem}[section]
\newtheorem{lemma}[theorem]{Lemma}

\newtheorem{definition}[theorem]{Definition}

\newtheorem{remark}[theorem]{Remark}

\newcounter{figures}[section]

\def\bE{{\mathbb E}}
\def\bN{{\mathbb N}}
\def\bF{{\mathbb F}}

\def\bR{{\mathbb R}}

\def\bX{{\mathbb X}}

\def\cC{\mathcal{C}}

\def\cK{\mathcal{K}}

\def\supp{{\rm{supp}\, }}

\def\ONE{{\mathbbm 1}}

\def\R{\mathbb{R}}

\begin{document}

\title[Minimax density estimation]
{Minimax density estimation on Sobolev spaces with dominating mixed smoothness}

\author{Galatia Cleanthous}
\address{School of Mathematics and Statistics,
Newcastle University, Newcastle Upon Tyne, NE1 7RU, UK}
\email{galatia.cleanthous@newcastle.ac.uk}

\author{Athanasios G. Georgiadis}
\address{Department of Mathematics and Statistics,
University of Cyprus, 1678 Nicosia, Cyprus}
\email{gathana@ucy.ac.cy}

\author{Emilio Porcu}
\address{School of Mathematics and Statistics,
Newcastle University, Chair of Spatial Analytics Methods Center}
\email{georgepolya01@gmail.com}



\begin{abstract}
We study minimax density estimation on the product space $\mathbb{R}^{d_1}\times\mathbb{R}^{d_2}$. We consider $L^p$-risk for probability density functions defined over regularity spaces that allow for different level of smoothness in each of the variables. Precisely, we study probabilities on Sobolev spaces with dominating mixed-smoothness. We provide the rate of convergence that is optimal even for the classical  Sobolev spaces.
\end{abstract}

\date{June 14, 2019}

\maketitle

\section{Introduction}\label{Introduction}
\setcounter{equation}{0}

\subsection{Context} 

\def\bX{\boldsymbol{X}}
Let $X$ be a random variable having a probability distribution that is absolutely continuous with respect to the Lebesgue measure, and with an unknown probability density function (denoted pdf throughout), $f$. A classical statistical problem is to estimate $f$, given a random independent identically distributed (iid) sample $\bX=(X_1,\dots,X_n)$ from $X$.

The minimax approach is a popular method for providing solutions to such a problem. We briefly sketch the idea behind the method: we call estimator of $f$, and denote it $\hat{f}$, a measurable function of the random vector $\bX$. We define the $L^p$-risk as $\mathbb{E}_f \|\hat{f}-f\|_p^p$, where $\mathbb{E}_f$ denotes expectation with respect to the probability measure $\mathbb{P}_f$ associated with $\bX$.

We assume that $f$ belongs to a {functional class} $\mathbb{F}$, and for a given estimator $\hat{f}$ the maximum risk is defined as the quantity
$$\sup_{f\in\mathbb{F}}\mathbb{E}_f \|\hat{f}-f\|_p^p.$$
Minimax approach is based on finding a rate optimal estimator $\mathbf{\hat{f}}$ such that
$$\sup_{f\in\mathbb{F}}\mathbb{E}_f \|\mathbf{\hat{f}}-f\|_p^p\sim\inf_{\tilde{f}}\sup_{f\in\mathbb{F}}\mathbb{E}_f \|\tilde{f}-f\|_p^p,$$
where the infimum is taken over all the possible estimators $\tilde{f}$. Then, $\mathbf{\hat{f}}$ is called the ``minimax estimator" with corresponding accuracy determined through the equation above. 

\subsection{Literature Review}

Minimax approaches have been popular within the non-parametric statistics research community  for many decades. The seminal paper by Bretagnolle and Huber \cite{BH} addressed the density estimation problem for pdf's having some regularity properties with respect to a given  Sobolev space $W^s_p (\bR)$. Tsybakov's book \cite{Tsybakov}, as well as the lecture notes \cite{HGPT} by H\"{a}rdle {\em et al.} provide a thorough introduction to density estimation within the nonparametric framework. A wealth of relevant contributions has been available in the next twenty years, and the reader is referred to  \cite{BKMP,CGKPP,DG,DL,DL2,DL3,DJKP,EY,EY2,GL2,GL3,Gol,HI,IK,KeP,KePT,KeLP,Rig,RigT,ST}. \cite{GL} provides a detailed historical overview of the research in this area up to the present decade. 

Apparently, solutions to the minimax problem rely strongy on  the function space $\mathbb{F}$ where pdf's are defined. It is customarily assumed that the function $f$ is sufficiently regular, {\em i.e.} it belongs to a given smoothness space. Prominent examples of smoothness spaces are  Sobolev, H\"{o}lder, Nikol'skij and Besov spaces. 
 
 Moreover, the function $f$ might be defined over multi-dimensional spaces. For instance, the $d$-dimensional Euclidean space, the cube, or the $(d-1)$-dimensional sphere embedded in $\mathbb{R}^d$. The paper by \cite{BKMP}
studies a similar problem for $f$ being defined on the sphere. More recently, \cite{CGKPP} considered this problem when $f$ is defined over manifolds or over more general metric spaces. 
\subsection{Our Contribution} 

We consider pdf's defined over the $d$-dimensional Euclidean space, $\bR^d$. For $d_1$, $d_2$ positive integers such that $d_1+d_2=d$, we consider functions $f$ having different orders of regularity over the two {\em directions} $\R^{d_1}$ and $\R^{d_2}$, respectively.  In nonparametric statistics, and in mathematical analysis, this case of different regularities over different directions is known as  \textit{mixed smoothness}. In particular, the study of spaces with mixed-smoothness goes back to the early '60s, with the fundamental contributions of the Russian school (see, for example, \cite{A,BIKLN,BIN,LN,nik,nik2,ScT}). Also the survey \cite{Sc} contains a full historical overview of the study of such spaces.

The problem of estimating a pdf on spaces with mixed smoothness, attracted significant attention inside the statistical community. The papers \cite{GL,GL2,GL3,HI,IK,KeLP,KLP2} challenge the problem of mixed smoothness under the name of \textit{anisotropic spaces}. Mixed smoothness received attention in spatial and space-time statistics as well: a covariance function might have different orders of differentiability over different directions in space. Also, normally spatial and temporal smoothness are different, as noted  by \cite{stein1} and subsequently by \cite{stein2}.

{\vspace{0.2cm}}

This paper introduces Sobolev spaces with \textit{dominating mixed smoothness}. To the knowledge of the authors, such a spaces have received a very limited attention in the statistical community. So much so, that we could not find any paper where an explicit use of these spaces has been advocated. 

{\vspace{0.2cm}}

Let $d_1,d_2$ and $s_1,s_2$ be positive integers. We study minimax density estimation for functions $f$ defined over products of Euclidean spaces $\bR^{d_1}\times\bR^{d_2}$, lying in Sobolev spaces with dominating mixed smoothness, denoted by $W^{(s_1,s_2)}_p (\bR^{d_1}\times\bR^{d_2})$. 

Surprisingly, it will turn out that these spaces support the mixed smoothness and simultaneously allow to provide bounds which are rate-optimal even for the ordinary (umnixed) Sobolev spaces.

This in turn implies that, for instance, for a pdf, $f$, that is not regular over one direction, {\em sharp} estimation can be achieved when  $f$ is smoother in the other direction. We show that higher smoothness in a given direction allows to compensate lower smoothness in the other direction. 

%
%
The remainder of the paper is as follows: In Section \ref{Sec2} we present the necessary analytical and statistical background, including the precise definition of the spaces we will work on. In Section \ref{Sec3} we provide some motivation and we state our main results. Theorems \ref{TH:UBpbig}, \ref{TH:UBpsmall2} and \ref{TH:UBpsmall1} include the upper bounds and Theorems \ref{th:lowerbound}, \ref{th:lowerboundnoncompact} the corresponding lower bounds. The upper bounds are obtained by a broad class of kernel density estimators which attain the optimal rate of convergence. Section \ref{Sec3} is accompanied with remarks and comparison of our results with classical and anisotropic spaces. Sections \ref{Sec4} and \ref{Prooflower} contain the proofs of our Theorems (upper and lower bounds respectively). For reasons of completeness in Section \ref{Sec6} we present some remarks on kernels and Sobolev spaces.

Let us summarize the contributions of our study:

$(\alpha)$ Theorems \ref{TH:UBpbig} and \ref{TH:UBpsmall2} deal with \textit{kernel} density estimators. An upper bound for the $L^p$-risk for pdf's on Sobolev spaces with dominating mixed-smoothness $W^{(s_1,s_2)}_p (\bR^{d_1}\times\bR^{d_2})$ is provided. 

$(\beta)$ Theorem \ref{th:lowerbound} provides the corresponding lower bound, concluding that the estimation is \textit{minimax}.

$(\gamma)$ For $1\le p<2$ and for a class of pdf's that may be non-compactly supported, an upper bound  is provided by Theorem \ref{TH:UBpsmall1}.

$(\delta)$ The precise behaviour of the lower bound for non-compactly supported pdf's when $1\le p<2$ is obtained in Theorem \ref{th:lowerboundnoncompact}.

$(\varepsilon)$ The minimax density estimation problem on classical Sobolev spaces $W^s_p(\bR^d)$ is presented in Theorems \ref{th:isotropic} and \ref{th:isotropic2}.

\section{Background} \label{Sec2}
This material is largely expository and provides the necessary ingredients to understand the theoretical results provided in the paper, as well as their proofs. 
\subsection{Analysis background} 
We start with some analysis concepts. 

\noindent {\em Minkowski's inequality.} We recall the generalized Minkowski's inequality: Let $1\le p\le\infty$. For any function $f$ defined over some product space, we write $\|f(x,y)\|_{p}$ for the $L_p$-norm with respect to $y$ for a fixed $x$. We have 
\begin{equation}
\label{Mink}
\Bigg \|\int_{\bR^d}f(x,y)dx\Bigg \|_p\le \int_{\bR^d}\|f(x,y)\|_p dx.
\end{equation}
{\em Young's inequality.} Recall that the {convolution} of two functions $f,g:\bR^d\rightarrow\bR$, is defined as
\begin{equation}
\label{Conv}
(f\ast g)(x):=\int_{\bR^d}f(x-y)g(y)dy,\quad\text{for every}\;\;x\in\bR^d.
\end{equation}
When $f\in L^p$, $p\in[1,\infty]$ and $g\in L^1$, by Young's inequality their convolution $f\ast g$ belongs to $L^p$. More precisely,
\begin{equation}
\label{Young}
\|f\ast g\|_p\le\|g\|_1\|f\|_p.
\end{equation}

\noindent {\em Multi-variable Taylor's formula} Let us fix some well-known multi-variable notation first. We denote by $\bN$ and $\bN_0$ the sets of positive and non-negative integers respectively. Let $x\in\bR^d$ and $\alpha\in\bN_0^d$. We denote by $x^{\alpha}=x_1^{\alpha_1}\cdots x_d^{\alpha_d}$, $|\alpha|:=\alpha_1+\cdots+\alpha_d$, the length of multi-index $\alpha$ and by $\alpha!=\alpha_1 !\cdots\alpha_d !$, its factorial.

Our action will take place on product spaces $\bR^{d_1}\times\bR^{d_2}$. We use a version of {Taylor's theorem} for functions having different levels of regularity corresponding to 
$\bR^{d_1}$ and $\bR^{d_2}$.

Let $s_1,s_2\in\bN$. A function $f:\bR^{d_1}\times\bR^{d_2}\rightarrow\bR$ belongs to the class $\cC^{(s_1,s_2)}=\cC^{(s_1,s_2)}(\bR^{d_1}\times\bR^{d_2})$ when the derivative $\partial^{\alpha}f=\partial_2^{\alpha_2}\partial_1^{\alpha_1}f$ is continuous, for every multi-index $\alpha=(\alpha_1,\alpha_2)\in\bN_0^{d_1}\times\bN_0^{d_2}$ such that $|\alpha_1|\le s_1$ and $|\alpha_2|\le s_2$.

Let now $s_1,s_2\in\bN$ and $f\in\cC^{(s_1,s_2)}$. Then, for every $(z,w),(x,y)\in\bR^{d_1}\times\bR^{d_2}$ we apply Taylor's formula on $\bR^{d_1}$ and subsequently on $\bR^{d_2}$ to get:

\begin{align}
\label{Taylor}
f(z,w)&=\sum_{|\alpha_1|<s_1,|\alpha_2|<s_2}\frac{\partial^{\alpha}f(x,y)}{\alpha!}(z-x)^{\alpha_1}(w-y)^{\alpha_2}
\\
&+\sum_{|\alpha_1|<s_1,|\alpha_2|=s_2}\frac{|\alpha_2|}{\alpha!}(z-x)^{\alpha_1}(w-y)^{\alpha_2}\int_{0}^{1}(1-t_2)^{|\alpha_2|-1}\partial^{\alpha}f(x,y+t_2(w-y))dt_2
\nonumber
\\
&+\sum_{|\alpha_1|=s_1,|\alpha_2|<s_2}\frac{|\alpha_1|}{\alpha!}(z-x)^{\alpha_1}(w-y)^{\alpha_2}\int_{0}^{1}(1-t_1)^{|\alpha_1|-1}\partial^{\alpha}f(x+t_1(z-x),y)dt_1
\nonumber
\\
&+\sum_{|\alpha_1|=s_1,|\alpha_2|=s_2}\frac{|\alpha_1||\alpha_2|}{\alpha!}(z-x)^{\alpha_1}(w-y)^{\alpha_2}\times\nonumber
\\
&\times\int_{0}^{1}\int_{0}^{1}(1-t_1)^{|\alpha_1|-1}(1-t_2)^{|\alpha_2|-1}\partial^{\alpha}f(x+t_1(z-x),y+t_2(w-y))dt_1 dt_2.
\nonumber
\end{align}

Note that this expression can be generalized to more general product spaces at the expense of very complicated notation. To avoid mathematical obfuscation, we work on the product of two spaces only, albeit our methods can be extended to the case of products of arbitrary many spaces, as illustrated in Section \ref{S:nproduct}.


Let us now recall the definition of Sobolev regularity spaces:

\begin{definition}\label{D:SS}
Let $s\in\mathbb{N}$, $1\le p<\infty$ and let $f$ be a function on $L^p (\bR^d)$. Then, $f$ belongs to the Sobolev space $W^s_p=W^s_p(\bR^d)$, when
\begin{equation}\label{Sob1}
\|f\|_{W^s_p}:=\|f\|_{W^s_p(\bR^d)}:=\sum_{|\alpha|\le s}\|\partial^{\alpha}f\|_{p}<\infty.
\end{equation}
\end{definition}

We shall deal with the following smoothness spaces with mixed smoothness on the product domain $\bR^{d_1}\times\bR^{d_2}$, which are the so called \textit{Sobolev spaces with dominating mixed smoothness}: 
\begin{definition}\label{D:LS}
Let $s_1,s_2\in\mathbb{N}$, $1\le p<\infty$ and let $f$ be a function on $\bR^{d_1}\times\bR^{d_2}$. Then, $f$ belongs to  Sobolev space $W^{(s_1,s_2)}_p=W^{(s_1,s_2)}_p(\bR^{d_1}\times\bR^{d_2})$, when
\begin{equation}\label{Sob2}\|f\|_{W^{(s_1,s_2)}_p}:=\|f\|_{W^{(s_1,s_2)}_p(\bR^{d_1}\times\bR^{d_2})}:=\sum_{|\alpha_1|\le s_1, |\alpha_2|\le s_2}\|\partial_2^{\alpha_2}\partial_1^{\alpha_1}f\|_{p}<\infty.
\end{equation}
Further we denote by $W^{(s_1,s_2)}_p(r)$ the closed ball of radius $r>0$ centered at the zero function, {\em i.e.}, $W^{(s_1,s_2)}_p(r)=\{f\in W^{(s_1,s_2)}_p:\;\|f\|_{W^{(s_1,s_2)}_p}\le r\}$.
\end{definition}
%

\begin{remark}\label{remark:inclusion} 
Some comments are in order: \\ -
$(\alpha)$ Let $s_1,s_2\in\bN$, $1\le p<\infty$ and $r>0$. The inequalities
\begin{equation}
\label{normsSobolev}
\|f\|_{W^{\min(s_1,s_2)}_p}\le\|f\|_{W^{(s_1,s_2)}_p}\le \|f\|_{W^{s_1+s_2}_p},
\end{equation}
imply the inclusion relations
\begin{equation}
\label{inclusionsSobolev}
W^{s_1+s_2}_p\subset W^{(s_1,s_2)}_p\subset W^{\min(s_1,s_2)}
\end{equation}
which apparently apply to the corresponding balls 
\begin{equation}
\label{ballsSobolev}
W^{s_1+s_2}_p(r)\subset W^{(s_1,s_2)}_p(r)\subset W^{\min(s_1,s_2)}(r).
\end{equation}
Thus, Sobolev spaces  with dominating mixed smoothness are embedded between classical Sobolev spaces of minimum smoothness $\min(s_1,s_2)$, and classical Sobolev spaces with a smoothness index being identically equal to $s_1+s_2$. 

$(\beta)$ Let $d_1,d_2$ and $s_1,s_2$ be positive integers and $f_i\in W^{s_i}_{p}(\bR^{d_i})$, $i=1,2$. Then their tensor product $f:=f_1\otimes f_2\in W^{(s_1,s_2)}_p(\bR^{d_1}\times\bR^{d_2})$.

$(\gamma)$ We are interested in probability density functions contained in balls of mixed smoothness Sobolev spaces. This fact implies some suitable restrictions on the radius of the ball where the pdf is defined. Precisely, when $p=1$, the restriction $r\ge\|f\|_{W^{(s_1,s_2)}_1}\ge\|f\|_1=1$ is needed for a pdf to be well defined over $W^{(s_1,s_2)}_p(r)$. Moreover for $r=1$, the only pdf's that belong to the ball $W^{(s_1,s_2)}_1(1)$ are (piecewise) constants. Thus, we shall avoid this case,  and will consider balls of radius $r>1$ when $p=1$.
\end{remark}

\subsection{Statistics background}
We now collect the statistical background material we need in our study. 

%
%

We consider a normalized kernel $K:\mathbb{R}^d\rightarrow \mathbb{R}$ with
$\int_{\mathbb{R}^d} K(y)dy=1$. The function
\begin{equation}\label{KdeR}
\hat{f_n}(x)=\frac{1}{nh^d}\sum\limits_{i=1}^{n} K\Big(\frac{x_i-x}{h}\Big), \quad x\in\mathbb{R}^d
\end{equation}
is called the \textit{kernel density estimator (kde) associated with the kernel K}. The parameter $h=h(n)$ is the \textit{bandwidth} of $\hat{f_n}$.

 The following classical inequalities will be used in the manuscript. \\
$(\alpha)$ \textit{Bernstein's inequality:} Let $Y_1,\dots,Y_n$ independent random variables such that $\mathbb{E}(Y_i)=0$, $\mathbb{E}\big(Y_i^2\big)\le \sigma^2$ and $|Y_i|\le M$, for every $i=1,\dots,n$. Then, for every $v>0$,
\begin{equation}\label{Bernstein}
\mathbb{P}\Bigg(\Big|\frac{1}{n}\sum_{i=1}^n Y_i\Big|\ge v\Bigg)\le 2\exp\Bigg(\frac{-nv^2}{2(\sigma^2+Mv/3)}\Bigg).
\end{equation}
$(\beta)$ \textit{Rosenthal's inequality:} Let $p\ge2$ and $Y_1,\dots,Y_n$ independent random variables such that $\mathbb{E}(Y_i)=0$ and $\mathbb{E}\big(|Y_i|^p\big)<\infty$ for every $i=1,\dots,n$. There exists a constant $c(p)>0$ such that
\begin{equation}\label{Rosenthal}
\mathbb{E}\Big(\Big|\sum_{i=1}^n Y_i\Big|^p\Big)\le c(p)\Big(\sum_{i=1}^n \mathbb{E}\big(|Y_i|^p\big)+\Big(\sum_{i=1}^n \mathbb{E}(Y_i^2)\Big)^{p/2}\Big).
\end{equation}
\begin{remark} 
When $0<p\leq 2$, by convexity we have
\begin{equation}\label{Rosenthal2}
\mathbb{E}\Big(\Big|\sum_{i=1}^n Y_i\Big|^p\Big)\le \Big(\sum_{i=1}^n \mathbb{E}\big(Y_i^2\big)\Big)^{p/2}.
\end{equation}
\end{remark}

\subsubsection{Minimax density estimation on classical Sobolev spaces}\label{Sub:Isotropic} 
Although this paper works under the framework of Sobolev spaces with dominating mixed smoothness, it will be useful to resort some properties of ordinary Sobolev spaces. 

Let $2\le p<\infty$ and $s\in\bN$. Arguments in  \cite{BH} show that the optimal rate for the minimax risk over the Sobolev space $W^s_p (\bR)$ is identically equal to
\begin{equation}
\label{optimalisotropic}
\inf_{\tilde{f}}\sup_{f\in W^s_p(\bR) (r)}\mathbb{E}_f \|\tilde{f}-f\|_p^p\sim n^{-ps/(2s+1)}.
\end{equation}
Apparently, the optimal rate depends on the index $s$ of regularity associated with the function $f$. Also, the approximation is improved when $s$ increases.  

To the knowledge of the authors, the result in Equation (\ref{optimalisotropic}) has not been extended to the case $W^{s}_p (\bR^d)$. Actually, the dimension $d$ has some implication on the rate of convergence: recently,  \cite{CGKPP} studied the density estimation problem on a class of metric spaces that include the $d$-dimensional Euclidean space $\bR^d$. We rephrase a result from \cite{CGKPP} to make it consistent with this exposition. \\
Let $d\in\bN$, $2\le p<\infty$ and $s\in\bN$. Then, for every $r>0$, the upper bound
\begin{equation}
\label{upperisotropic}
\sup_{f\in W^s_p(\bR^d) (r)}\mathbb{E}_f \|\hat{f}_n-f\|_p^p\le c n^{-ps/(2s+d)},
\end{equation}
applies for a broad class of kernel density estimators $\hat{f}_n$ as defined through Equation (\ref{KdeR}). Apparently, the ratio above  depends on the dimension $d$ of the space where the pdf is defined. This might be expected: see, for instance, 
\cite{DJKP,GL,GL2,GL3,HI,IK}. For reasons of completeness we do present the minimax theorem for classical (unmixed) Sobolev spaces on $\bR^d$ in our last Section.

\section{Results} \label{Sec3}

In this Section we present the motivation for our study, we state our results and we compare our paper with other contributions in the area.

\subsection{Motivation}\label{motivation}

The bound (\ref{upperisotropic}) in concert with the inclusions in Remark \ref{remark:inclusion}, suggest the following upper bound for the minimax risk:

Let $d_1,d_2\in\bN$, $2\le p<\infty$ and $s_1,s_2\in\bN$. Then it turns out that
\begin{equation}
\label{sureanisotropic}
\sup_{f\in W^{(s_1,s_2)}_p (r)}\mathbb{E}_f \|\tilde{f}-f\|_p^p\le\sup_{f\in W^{s_{\min}}_p (r)}\mathbb{E}_f \|\tilde{f}-f\|_p^p\sim n^{-ps_{min}/(2s_{min}+(d_1+d_2))},
\end{equation}
where $s_{\min}:=\min(s_1,s_2)$.

{\vspace{0.2cm}}

Our main point is that such an upper bound might be suboptimal. A clear evidence is provided by the case where the the function $f$ is much smoother in one direction with respect to the other one. Clearly in (\ref{sureanisotropic}) the bound depends only on the variable in which $f$ is less smooth and we do not gain anything from the ``good" variable. Hence the need for studying the problem from the perspective of Sobolev spaces with mixed smoothness, with the hope that we can somehow improve (\ref{sureanisotropic}) by involving the $s_{\max}:=\max(s_1,s_2)$.

Combination of results in Section \ref{Sub:Isotropic} with the inclusions appearing in Remark \ref{remark:inclusion} explains that 
the best possible bound that one should expect for the $L^p$-risk under study should be equal to  
\begin{equation}
\label{bestupperbound}
\frac{s_1+s_2}{2(s_1+s_2)+(d_1+d_2)}.
\end{equation}
The way we approach this best possible exponent will become apparent subsequently. 

\subsection{Kernel density estimators on spaces with mixed smoothness}

\begin{definition}\label{kernelsclass} Let $s_1,s_2\in\bN$. A kernel $K:\mathbb{R}^{d_1}\times\mathbb{R}^{d_2}\rightarrow \mathbb{R}$ belongs to the class $\cK(s_1,s_2)$ when

1. Markov property:
\begin{equation}\label{K:Markov}
\iint_{\mathbb{R}^{d_1}\times\mathbb{R}^{d_2}} K(u)du=1.
\end{equation}
2. $K$ has vanishing moments of any order ; $1\le|\alpha|<s_1+s_2$ with $|\alpha_i|\le s_i,\;i=1,2$.
\begin{equation}\label{K:moments}
\iint_{\bR^{d_1}\times\bR^{d_2}}u^{\alpha}K(u)du=0.
\end{equation}
3. The following integrals are finite:
\begin{equation}
\label{K:finite}
\iint_{\bR^{d_1}\times\bR^{d_2}}|u_1^{\alpha_1}||u_2^{\alpha_2}||K(u_1,u_2)|du_1du_2\le I_{(s_1,s_2)}<\infty,
\end{equation}
for $|\alpha_1|=s_1$ and $|\alpha_2|=s_2.$ \\
4. The kernel is bounded:
\begin{equation}
\label{K:bounded}
\|K\|_{\infty}:=\sup_{u\in\bR^{d_1}\times\bR^{d_2}}|K(u)|<\infty.
\end{equation}
\end{definition}
Let $n\in\mathbb{N}$ and $(X_1,Y_1),\dots,(X_n,Y_n)$ be iid random variables, with probability density function $f$. We extend the definition of kernel density estimation to this product space through 
\begin{equation}\label{KdeProduct}
\hat{f_n}(x,y)=\frac{1}{nh^{d_1}h^{d_2}}\sum\limits_{i=1}^{n} K\Big(\frac{x_i-x}{h},\frac{y_i-y}{h}\Big), \quad (x,y)\in\mathbb{R}^{d_1}\times\mathbb{R}^{d_2},
\end{equation}
where $0<h=h_n<1$ is the bandwidth.

The existence of kernels belonging to the class $\cK(s_1,s_2)$ will be discussed in Section \ref{S:Kernels}
Note that we can invoke (\ref{K:Markov}), (\ref{K:finite}), (\ref{K:bounded}) in concert with Riesz-Thorin's Theorem to show that
\begin{equation}
\label{K:Lq}
\|K\|_q<\infty,\quad \text{for every}\;\;q\in[1,\infty].
\end{equation}

\subsection{Upper bounds}
We start by considering pdf's defined over mixed-smoothness Sobolev balls in the $L^p$-norm, using kernel density estimators generated by kernels of the above class. Our first result provides an upper bound for the case $p\ge2$.

\begin{theorem}\label{TH:UBpbig}
Let $p\ge2$, $r>0$, $s_1,s_2\in\bN$ and let a kernel $K$ belong to the class $\cK(s_1,s_2)$. Let $\hat{f}_n$ be the corresponding kernel density estimator defined as in (\ref{KdeProduct}). Then, there exists a constant $c>0$ such that
\begin{equation}\label{upper1}
\sup_{f\in W^{(s_1,s_2)}_p (r)}\mathbb{E}\|\hat{f}_n-f\|_p^p\le cn^{-(s_1+s_2)p/(2(s_1+s_2)+d_1+d_2)}.
\end{equation}
\end{theorem}

\begin{remark}
Some comments are in order. 

$(\alpha)$ A close look at the proof of Theorem \ref{TH:UBpbig} provides a uniquely determined value for the constant $c>0$:
\begin{equation}
\label{constant}
c=2^{p-1}\Bigg [\Bigg (I_{(s_1,s_2)}\sum_{|\alpha_1|=s_1,|\alpha_2|=s_2}\|\partial^{\alpha}f\|_p\Bigg)^p+c(p)2^{p-2}\|K\|_{\infty}^{p-2}\|K\|_2^2+c(p)\|K\|_2^p\|f\|_{p/2}^{p/2}\Bigg ],
\end{equation}
where $I_{(s_1,s_2)}$ is the constant from (\ref{K:finite}), and $c(p)$ was determined at (\ref{Rosenthal}). Note that the smoothness norm that appears above is the norm of $f$ on the so-called homogeneous Sobolev space with dominating mixed smoothness.

$(\beta)$ The rate we succeed in (\ref{upper1}) is exactly (\ref{bestupperbound}); the optimal one for the classical Sobolev space $W^{s_1+s_2}_p(\bR^{d_1+d_2})$, even though we worked for the bigger Sobolev space $W^{(s_1,s_2)}_{p}(\bR^{d_1}\times\bR^{d_2})$ with dominating mixed smoothness. As a conclusion, we can see that $s_{\max}$ performs as an antidote for $s_{\min}$ (see section \ref{motivation}).

\end{remark}

The case $p<2$ needs a separate treatment. In particular, the case $p=1$ requires a bunch of additional technical assumptions (see for example \cite{BH}). In the following, we assume that  $f$ is compactly supported.

Let $r,R>0$. We denote by $W^{(s_1,s_2)}_p (r,R)$ the set of all pdf's $f\in W^{(s_1,s_2)}_p (r)$ such that
$\supp f\subset \{(x,y)\in\bR^{d_1}\times\bR^{d_2}: |(x,y)-(x_0,y_0)|\le R\}$, for some $(x_0,y_0)\in\bR^{d_1}\times\bR^{d_2}$. We are now able to state the following result. 

\begin{theorem}
\label{TH:UBpsmall2}
Let $1\le p<2$, $r,R>0$, $s_1,s_2\in\bN$ and let $K\in\cK(s_1,s_2)$ be compactly supported. Let $\hat{f}_n$ be the corresponding kernel density estimator as defined in (\ref{KdeProduct}).
Then, there exists a constant $c>0$ such that
\begin{equation}
\sup_{f\in W^{(s_1,s_2)}_p (r,R)}\mathbb{E}\|\hat{f}_n-f\|_p^p\le cn^{-(s_1+s_2)p/(2(s_1+s_2)+d_1+d_2)}.
\end{equation}
\end{theorem}

\subsection{Minimax density estimation} So far, we have provided a way to construct estimators $\hat{f}_n$ with a maximum risk on a ball of a mixed smoothness Sobolev space that is bounded from above by a certain rate. Precisely,
\begin{equation}\label{upperrecall}
\sup_{f\in W^{(s_1,s_2)}_p (r)}\mathbb{E}\|\hat{f}_n-f\|_p^p\le cn^{-pS/(2S+D)},
\end{equation}
where for brevity we set $S:=s_1+s_2$ and $D:=d_1+d_2$.

{\vspace{0.3cm}}

The inspection for a lower bound starts by seeking for a constant $c>0$ such that
\begin{equation}\label{lowerboundintro}
\inf_{\tilde{f}}\sup_{f\in W^{(s_1,s_2)}_p (r)}\mathbb{E}\|\tilde{f}-f\|_p^p\ge cn^{-pS/(2S+D)},
\end{equation}
where the infimum is taken over all possible estimators $\tilde{f}$ and for $n$ sufficiently large.

{\vspace{0.3cm}}

\begin{remark}
There are some technical restrictions that appear when one works with pdf's defined over subspaces of the space $L^1$. As mentioned earlier, for $p=1$ the ball $W^{(s_1,s_2)}_1 (r)$ contains well-defined pdf's only when $r\ge1$. Further, for $r=1$ the only pdf`s that are well-defined have vanishing derivatives of all orders. This is the trivial case of (piecewise) constant pdf's that we may avoid in our study. So in the case when $p=1$ the radius $r$ will be always assumed to be greater than $1$. To unifying notation, we set
\begin{equation}\label{r*}
r_{*}:=
\begin{cases}
\;r-1\;,\;\;p=1\\
\;r,\;\;p>1.
\end{cases}
\end{equation}
\end{remark}
We prove the following lower bound.
\begin{theorem}
\label{th:lowerbound}
(i) Let $p\ge2$, $s_1,s_2\in\bN$ and $r>0$. 

Then, there exists a constant $c>0$ such that
\begin{equation}\label{lowerbound}
\liminf\limits_{n\rightarrow\infty}\Bigg(\frac{r^{D/S}}{n}\Bigg)^{-\frac{pS}{2S+D}}\inf_{\tilde{f}}\sup_{f\in W^{(s_1,s_2)}_p (r)}\mathbb{E}\|\tilde{f}-f\|_p^p\ge c.
\end{equation}

(ii) Let $1\le p<2$, $s_1,s_2\in\bN$, $r_{*}$ as in relation (\ref{r*}) and $R>0$ large enough. 

Then, there exists a constant $c>0$ such that
\begin{equation}\label{lowerboundsmallp}
\liminf\limits_{n\rightarrow\infty}\Bigg(\frac{r_{*}^{D/S}}{n}\Bigg)^{-\frac{Sp}{2S+D}}\inf_{\tilde{f}}\sup_{f\in W^{(s_1,s_2)}_p (r,R)}\mathbb{E}\|\tilde{f}-f\|_p^p\ge c.
\end{equation}
\end{theorem}
Note that the infimum is taken over all possible estimators $\tilde{f}$ and the constant $c$ is independent of $r$ and $r_*$ respectively.
The proof of Theorem \ref{th:lowerbound} is highly inspired by \cite{GL} and is provided in Section \ref{Prooflower}.

\subsection{Further results for non-compactly supported pdf's for the range \boldmath $1\le p<2$.}\label{ncpsmall}
This section shows that, when $1 \le p <2$, the assumption of compact support in Theorem \ref{TH:UBpsmall2} can be eluded if additional technicalities are assumed. A suggestion comes from Kerkyacharian and Picard in \cite{KeP}, who assume a pdf to be dominated by a radial and radially dicreasing $L^{p/2}$ bounded function. A similar result is provided here. Some additional definitions and notations are needed. 
\begin{definition}\label{D:classomega}
Let $\omega:\bR^{d_1}\times\bR^{d_2}\rightarrow[0,\infty)$ be a radial function that is radially non-increasing. Let $\omega$ belong to $L^{p/2}\cap L^{\infty}$. We denote by $\tilde{W}^{(s_1,s_2)}_{p}(r)$ the space of all pdf's that belong to the Sobolev ball $W^{(s_1,s_2)}_{p}(r)$ and that additionally satisfy the domination $f(x-x_0,y-y_0)\le \omega(x,y)$ for some fixed $(x_0,y_0)\in\bR^{d_1}\times\bR^{d_2}$ and for some $\omega$ being defined as above.
\end{definition}
\noindent We are now able to state our result.
\begin{theorem}
\label{TH:UBpsmall1}
Let $p,r,s_1,s_2$ and $K$ be as in Theorem \ref{TH:UBpsmall2}. \\
Then, there exists a constant $c>0$ such that
\begin{equation}
\sup_{f\in \tilde{W}^{(s_1,s_2)}_p (r)}\mathbb{E}\|\hat{f}_n-f\|_p^p\le cn^{-(s_1+s_2)p/(2(s_1+s_2)+d_1+d_2)}.
\end{equation}
\end{theorem}
Some comments are in order.
Theorems \ref{TH:UBpbig} and \ref{th:lowerbound} show that, for $p>2$, the rate $n^{-S/(2S+D)}$ is minimax for the class $W^{(s_1,s_2)}_{p} (r)$. However, for $1\le p<2$ we achieve minimax estimation only for compactly supported pdf's, as it can be verified by Theorems \ref{TH:UBpsmall2} and \ref{th:lowerbound}. \\
We focus again on how to elude the assumption of compact support while trying to improve the lower bound  $1\le p<2$. The result following subsequently illustrates our findings. 
\begin{theorem}
\label{th:lowerboundnoncompact}
Let $1\le p<2$, $s_1,s_2\in\bN$ and $r_{*}$ as in (\ref{r*}).  \\
Then, there exists a constant $c>0$ such that
\begin{equation}\label{lowerboundsmallpnoncompact}
\liminf\limits_{n\rightarrow\infty}\Big(\frac{r_{*}^{D/S}}{n}\Big)^{-\frac{S(p-1)p}{Sp+D(p-1)}}\inf_{\tilde{f}}\sup_{f\in W^{(s_1,s_2)}_p (r)}\mathbb{E}\|\tilde{f}-f\|_p^p\ge c,
\end{equation}
where the infimum is taken over all possible estimators $\tilde{f}$ and the constant $c$ is independent of $r$.
\end{theorem}
The proof is technical and deferred to Section \ref{Prooflower}.


\subsection{Comparison with other smoothness spaces}\label{Remarksection}

\subsubsection{Comparison with classical Sobolev spaces} 

The discussion in Section \ref{motivation} shows that the rate to be expected should lie within two extremes being, respectively, 
$$ \frac{s_{\min}}{2s_{min}+d_1+d_2} \qquad \text{and} \qquad \frac{s_1+s_2}{2(s_1+s_2)+(d_1+d_2)} ,$$
and apparently the upper limit of this interval is the best possible. Therefore, here we are able to {match the best possible bound for a broader class of pdf's} belonging to $W^{(s_1,s_2)}_p (\bR^{d_1}\times\bR^{d_2})$ rather than $W^{s_1+s_2}_p(\bR^{d_1+d_2})$.

 Figure \ref{Nasos_figure} illustrates our framework: we depict 
 the derivatives that need to be integrable on $L^p$ for a function $f$ belonging either to the Sobolev space with dominating mixed smoothness, or to the corresponding classical Sobolev space. Here, we set $d_1=d_2=1$ and $s_1=1$, $s_2=4$, so that $s_1+s_2=5$. In particular, using the results available for the classical case, when working with the class $W^{(4,1)}_{p}(\bR^2)$ the rate would be $n^{-1/4}$, which is clearly improved by our rate, that in this case would be  $n^{-5/12}$.

\begin{figure} \label{Nasos_figure}
\centering
\hspace{-2cm}
\begin{minipage}{.47\textwidth}
\begin{tikzpicture}
\draw[gray, thick] (0,1) -- (4,1);
\draw[gray, thick](4,1)--(4,0);
\filldraw (0,0) circle (0.08cm);
\filldraw (0,1) circle (0.08cm)node[left] {1} ;

\filldraw (1,0) circle (0.08cm);
\filldraw (1,1) circle (0.08cm);

\filldraw (2,0) circle (0.08cm);
\filldraw (2,1) circle (0.08cm);

\filldraw (3,0) circle (0.08cm);
\filldraw (3,1) circle (0.08cm);

\filldraw (4,0) circle (0.08cm) node[below] {4};
\filldraw (4,1) circle (0.08cm);


\draw[->] (0,0) -- (4.5,0) node[below] {$\frac{\partial}{\partial x_1}$};
\draw[->] (0,0) -- (0,4.5) node[below left] {$\frac{\partial}{\partial x_2}$};

\end{tikzpicture}
\end{minipage}
\hspace{0.5cm}
\begin{minipage}{.4\textwidth}
\begin{tikzpicture}
\draw[gray, thick] (0,5) -- (5,0);


\filldraw (0,0) circle (0.08cm);
\filldraw (0,1) circle (0.08cm) node[left] {1} ;
\filldraw (0,2) circle (0.08cm);
\filldraw (0,3) circle (0.08cm);
\filldraw (0,4) circle (0.08cm);
\filldraw (0,5) circle (0.08cm) node[left] {5};

\filldraw (1,0) circle (0.08cm);
\filldraw (1,1) circle (0.08cm);
\filldraw (1,2) circle (0.08cm);
\filldraw (1,3) circle (0.08cm);
\filldraw (1,4) circle (0.08cm);

\filldraw (2,0) circle (0.08cm);
\filldraw (2,1) circle (0.08cm);
\filldraw (2,2) circle (0.08cm);
\filldraw (2,3) circle (0.08cm);

\filldraw (3,0) circle (0.08cm);
\filldraw (3,1) circle (0.08cm);
\filldraw (3,2) circle (0.08cm);

\filldraw (4,0) circle (0.08cm) node[below] {4} ;
\filldraw (4,1) circle (0.08cm);

\filldraw (5,0) circle (0.08cm) node[below] {5} ;


\draw[-] (0,1) -- (4,1);
\draw[-] (4,0) -- (4,1);

\draw[->] (0,0) -- (6.5,0) node[below] {$\frac{\partial}{\partial x_1}$};
\draw[->] (0,0) -- (0,6.5) node[below left] {$\frac{\partial}{\partial x_2}$};
\end{tikzpicture}
\end{minipage}
\caption{Left: $W^{(4,1)}_p(\bR\times\bR)$.  Right: $W^{5}_p(\bR^2)$.}
\end{figure}

\subsubsection{Comparison with anisotropic Sobolev spaces} 
Spaces with mixed smoothness are used extensively the last years by statistics community. In \cite{GL,GL3} Goldenshluger and Lepski study density estimation in \textit{anisotropic Nikol'skij spaces} $B^{(s_1,s_2)}_{p,\infty}$ (being a special case of Besov spaces). Let us compare their results expressed for \textit{anisotropic Sobolev spaces} (for consistency and simplicity), with ours.

\begin{definition}\label{D:AS}
Let $s_1,s_2\in\mathbb{N}$, $1\le p<\infty$ and let $f$ be a function on $\bR^{d_1}\times\bR^{d_2}$. Then, $f$ belongs to anisotropic Sobolev space $W^{(s_1,s_2)}_{p;\text{{\rm aniso}}}=W^{(s_1,s_2)}_{p;\text{{\rm aniso}}}(\bR^{d_1}\times\bR^{d_2})$, when
\begin{equation}\label{Sobaniso}\|f\|_{W^{(s_1,s_2)}_{p;\text{aniso}}}:=\sum_{0\le\frac{|\alpha_1|}{s_1}+\frac{|\alpha_2|}{s_2}\le1}\|\partial_2^{\alpha_2}\partial_1^{\alpha_1}f\|_{p}<\infty.
\end{equation}
\end{definition}
The following inclusion relation sheds a light on the comparison provided in this section \cite{ST}:
\begin{equation}
\label{incandmcl}
W^{s_1+s_2}_{p}(\bR^{d_1+d_2})\subset W^{(s_1,s_2)}_{p}(\bR^{d_1}\times\bR^{d_2})\subset W^{(s_1,s_2)}_{p;\text{aniso}}(\bR^{d_1}\times\bR^{d_2})\subset W^{s_{\min}}_{p}(\bR^{d_1+d_2}).
\end{equation}
This inclusion is depicted by Figure \ref{Figure_2}, where we can see the derivatives that belong to $L^p$ for the three spaces when $d_1=d_2=1$, $s_1=4$ and $s_2=1$.

\begin{figure} \label{Figure_2}
\centering
\hspace{-2cm}

\begin{tikzpicture}
\draw[gray, thick] (0,5) -- (5,0);
\draw[gray, thick] (0,1) -- (4,0);


\filldraw (0,0) circle (0.08cm);
\filldraw (0,1) circle (0.08cm) node[left] {1} ;
\filldraw (0,2) circle (0.08cm);
\filldraw (0,3) circle (0.08cm);
\filldraw (0,4) circle (0.08cm);
\filldraw (0,5) circle (0.08cm) node[left] {5};

\filldraw (1,0) circle (0.08cm);
\filldraw (1,1) circle (0.08cm);
\filldraw (1,2) circle (0.08cm);
\filldraw (1,3) circle (0.08cm);
\filldraw (1,4) circle (0.08cm);

\filldraw (2,0) circle (0.08cm);
\filldraw (2,1) circle (0.08cm);
\filldraw (2,2) circle (0.08cm);
\filldraw (2,3) circle (0.08cm);

\filldraw (3,0) circle (0.08cm);
\filldraw (3,1) circle (0.08cm);
\filldraw (3,2) circle (0.08cm);

\filldraw (4,0) circle (0.08cm) node[below] {4} ;
\filldraw (4,1) circle (0.08cm);

\filldraw (5,0) circle (0.08cm) node[below] {5} ;


\draw[-] (0,1) -- (4,1);
\draw[-] (4,0) -- (4,1);

\draw[->] (0,0) -- (6.5,0) node[below] {$\frac{\partial}{\partial x_1}$};
\draw[->] (0,0) -- (0,6.5) node[below left] {$\frac{\partial}{\partial x_2}$};
\end{tikzpicture}

\caption{From left to right: $W^{(4,1)}_{p;\text{aniso}}(\bR\times\bR)$, $W^{(4,1)}_p(\bR\times\bR)$ and $W^{5}_p(\bR^2)$.}
\end{figure}

However, the corresponding bounds that have been achieved by \cite{GL}, apply to our case as well, and precisely will be of order
\begin{equation}
\label{boundlepski}
\frac{1}{2+\frac{d_1}{s_1}+\frac{d_2}{s_2}},
\end{equation}
which is strictly smaller than (\ref{bestupperbound}), since the anisotropic space is bigger than the one a space with dominating mixed smoothness. For the special case we mentioned before, $d_1=d_2=1$, $s_1=4$ and $s_2=1$, the rate coming from the bound (\ref{boundlepski}) equals $n^{-4/13}$, which is bigger than $n^{-5/12}$.

Therefore, the ratio provided in this paper slightly outperforms the one in \cite{GL}. This fact strengthen our choice and motivation to study spaces with \textit{dominating mixed smoothness}.

We now provide a conclusive comment. The better rate obtained for the case
$W^{(4,1)}_p(\bR\times\bR)$ with respect to the case $W^{(4,1)}_{p;aniso}(\bR\times\bR)$, is  justified by the following equivalent norms for the two spaces (see for example \cite{ST}):
\begin{equation}
\label{equivaniso}
\|f\|_{W^{(s_1,s_2)}_{p;aniso}}\sim \|f\|_p+\sum_{|\alpha_1|=s_1}\|\partial_1^{\alpha_1}f\|_{p}+\sum_{|\alpha_2|=s_2}\|\partial_2^{\alpha_2}f\|_{p},
\end{equation}

\begin{equation}
\label{equidms}
\|f\|_{W^{(s_1,s_2)}_{p}}\sim \|f\|_p+\sum_{|\alpha_1|=s_1}\|\partial_1^{\alpha_1}f\|_{p}+\sum_{|\alpha_2|=s_2}\|\partial_2^{\alpha_2}f\|_{p}+\sum_{|\alpha_1|=s_1,|\alpha_2|=s_2}\|\partial_2^{\alpha_2}\partial_1^{\alpha_1} f\|_{p}.
\end{equation}

This use of the mixed derivatives $\partial_2^{\alpha_2}\partial_1^{\alpha_1} f$ of order $s_1+s_2$ was exactly the reason why in $W^{(s_1,s_2)}_{p}(\bR^{d_1}\times\bR^{d_2})$, we obtain a better bound than in $W^{(s_1,s_2)}_{p;aniso}(\bR^{d_1}\times\bR^{d_2})$. The analogous situation for Besov spaces, is the extra use of mixed-differences (see \cite[P. 156]{Sc}).

\section{Proofs of the upper bounds} \label{Sec4}

In this section we present the proofs of Theorems \ref{TH:UBpbig}, \ref{TH:UBpsmall1} and \ref{TH:UBpsmall2}. We start by noting that the risk $\mathbb{E}_f\|\hat{f}_n-f\|_p^p$ can be decomposed into two terms:
\begin{align}\label{Tr}
\mathbb{E}_f\big\|\hat{f}_n-f\big\|_p^p\leq 2^{p-1} \Big(\big\|\mathbb{E}_f[\hat{f}_n]-f\big\|_p^p+\mathbb{E}\big\|\hat{f}_n-\mathbb{E}_f[\hat{f}_n]\big\|^p_p\Big)
&=:2^{p-1}(B+S),
\end{align}
where $B$ and $S$ are the \textit{bias} and \textit{stochastic} terms respectively. We study these two terms separately. Note that the term $B$ determines the regularity of the pdf, $f$. It is notorious that there is a trade off between and $B$ and $S$, and the way such a trade off is balanced is through a proper choice of the bandwidth. We also note that a different approach is needed for the stochastic term, $S$, depending on whether  $p\ge2$ or $p<2$. This fact is justified by technical arguments coming from both analysis and statistics. Specifically,  when $p\ge2$, by interpolation Riesz-Thorin's Theorem, a pdf that lies in $W^{(s_1,s_2)}_{p}$ belongs to $L^{p/2}$ as well, which is not true for $p<2$. On the other hand, for $p\ge2$, we can invoke Rosenthal's inequality (\ref{Rosenthal}). 

\subsection{Estimation of Bias}
Let $b(x,y):=\mathbb{E}_f\big[\hat{f_n}(x,y)\big]-f(x,y)$. Then, by (\ref{K:Markov}) we get
\begin{equation}\label{pr1.2}
b(x,y)=\frac{1}{h^{d_1+d_2}}\iint_{\mathbb{R}^{d_1}\times\mathbb{R}^{d_2}} K\Big(\frac{z-x}{h},\frac{w-y}{h}\Big)(f(z,w)-f(x,y))dzdw.
\end{equation}

By Taylor's formula (\ref{Taylor}) in concert with the assumptions (\ref{K:moments}) and some change of variables, we derive

\begin{align}\label{pr1.3}
b(x,y)&=\sum_{|\alpha_1|=s_1,|\alpha_2|=s_2}\frac{|\alpha_1||\alpha_2|}{\alpha!}\frac{1}{h^{d_1+d_2}}\iint_{\mathbb{R}^{d_1}\times\mathbb{R}^{d_2}} K\Big(\frac{z-x}{h},\frac{w-y}{h}\Big)(z-x)^{\alpha_1}(w-y)^{\alpha_2}\times
\nonumber
\\
&\times\int_{0}^{1}\int_{0}^{1}(1-t_1)^{|\alpha_1|-1}(1-t_2)^{|\alpha_2|-1}\partial^{\alpha}f(x+t_1(z-x),y+t_2(w-y))dt_1 dt_2dzdw
\nonumber
\\
&=\sum_{|\alpha_1|=s_1,|\alpha_2|=s_2}\frac{|\alpha_1||\alpha_2|}{\alpha!}h^{|\alpha_1|+|\alpha_2|}\iint_{\mathbb{R}^{d_1}\times\mathbb{R}^{d_2}} K(z,w)z^{\alpha_1}w^{\alpha_2}\times
\nonumber
\\
&\times\int_{0}^{1}\int_{0}^{1}(1-t_1)^{|\alpha_1|-1}(1-t_2)^{|\alpha_2|-1}\partial^{\alpha}f(x+t_1hz,y+t_2hw)dt_1 dt_2dzdw.
\nonumber
\end{align}

We can now apply the triangle and (generalized) Minkowski's inequality (\ref{Mink}) in concert with the assumption (\ref{K:finite}) to obtain

\begin{align*}
\|b(x,y)\|_p&\le\sum_{|\alpha_1|=s_1,|\alpha_2|=s_2}\frac{|\alpha_1|\alpha_2|}{\alpha!}h^{s_1+s_2}\iint_{\mathbb{R}^{d_1}\times\mathbb{R}^{d_2}} |K(z,w)||z^{\alpha_1}||w^{\alpha_2}|\times
\nonumber
\\
&\times\int_{0}^{1}\int_{0}^{1}(1-t_1)^{|\alpha_1|-1}(1-t_2)^{|\alpha_2|-1}\|\partial^{\alpha}f(x+t_1hz,y+t_2hw)\|_p dt_1 dt_2dzdw
\nonumber
\\
&=\sum_{|\alpha_1|=s_1,|\alpha_2|=s_2}\frac{|\alpha_1|\alpha_2|}{\alpha!}h^{s_1+s_2}\iint_{\mathbb{R}^{d_1}\times\mathbb{R}^{d_2}} |K(z,w)||z^{\alpha_1}||w^{\alpha_2}|dzdw\times
\nonumber
\\
&\times\|\partial^{\alpha}f\|_p\int_{0}^{1}(1-t_1)^{|\alpha_1|-1}dt_1\int_{0}^{1}(1-t_2)^{|\alpha_2|-1}  dt_2
\nonumber
\\
&\le I_{(s_1,s_2)}h^{s_1+s_2}\sum_{|\alpha_1|=s_1,|\alpha_2|=s_2}\|\partial^{\alpha}f\|_p.
\end{align*}

This allows to conclude that the bias, $B$, is bounded by the quantity
\begin{equation}
\label{Biasupperbound}
B\le I_{(s_1,s_2)}^p h^{(s_1+s_2)p}\Big(\sum_{|\alpha_1|=s_1,|\alpha_2|=s_2}\|\partial^{\alpha}f\|_p\Big)^p\le ch^{(s_1+s_2)p}\|f\|_{W^{(s_1,s_2)}_p}^p.
\end{equation}

\subsection{Estimation of Stochastic term} 
We are now ready to estimate the stochastic term. Let us fix the real number $h$ such that
\begin{equation}
\label{h:choise}
h= n^{-1/(2(s_1+s_2)+(d_1+d_2))}.
\end{equation}
We set the random variables for every $i=1,\dots,n$.
\begin{equation*}
\eta_i(x,y):=K\Big(\frac{x_i-x}{h},\frac{y_i-y}{h}\Big)-\mathbb{E}_f\Big[K\Big(\frac{x_i-x}{h},\frac{y_i-y}{h}\Big)\Big],\quad (x,y)\in \bR^{d_1}\times\bR^{d_2}.
\end{equation*}

The random variables $\eta_1(x,y),\dots,\eta_n(x,y)$ are iid with $\mathbb{E}[\eta_i(x,y)]=0,$ for every $i=1,\dots,n$ and their variance is bounded by
\begin{eqnarray}\label{eYi22}
\mathbb{E}_f[\eta_i^2(x,y)]&\leq& \mathbb{E}_f\Big[K^2\Big(\frac{x_i-x}{h},\frac{y_i-y}{h}\Big)\Big] 
\\ \nonumber
&=&\iint_{\bR^{d_1}\times\bR^{d_2}} K^2\Big(\frac{z-x}{h},\frac{w-y}{h}\Big) f(z,w) dzdw.
\end{eqnarray}

\subsubsection{Proof of Theorem \ref{TH:UBpbig}} Let $2\leq p < \infty$. We observe that $\mathbb{E}_f\big(\eta_i(x,y)\big)=0$. By Fubini-Tonelli Theorem and Rosenthal's inequality, we derive
\begin{align}\label{S}\nonumber
S&=\mathbb{E}_f\big\|\hat{f}_n-\bE_f[\hat{f}_n]\big\|_p^p=\mathbb{E}_f\Big\|\frac{1}{n}\frac{1}{h^{d_1+d_2}}\sum\limits_{i=1}^{n}\eta_i(x,y)\Big\|_p^p
\\ \nonumber
&=n^{-p}h^{-p(d_1+d_2)}\iint_{\bR^{d_1}\times\bR^{d_2}}\mathbb{E}_f\Big(\Big|\sum_{i=1}^n \eta_i(x,y)\Big|^p\Big) dxdy\\ \nonumber
&\leq c(p)n^{-p}h^{-p(d_1+d_2)}\iint_{\bR^{d_1}\times\bR^{d_2}}\sum_{i=1}^n\mathbb{E}_f\big(|\eta_i(x,y)|^p\big) dxdy\\ \nonumber
&+
c(p)n^{-p}h^{-p(d_1+d_2)}\iint_{\bR^{d_1}\times\bR^{d_2}}\Big(\sum_{i=1}^n\mathbb{E}_f\big(\eta_i^2(x,y)\big)\Big)^{p/2} dxdy\\ 
&=:c(p)(S_1+S_2).
\end{align}
We shall work with $S_1$ and $S_2$ separately.

{\vspace{0.3cm}}

\textit{Estimation of $S_1$.}  \\
First, note that
\begin{equation}
\label{esup}
|\eta_i (x,y)|=\Big|K\Big(\frac{x_i-x}{h},\frac{y_i-y}{h}\Big)-\mathbb{E}_f\Big[K\Big(\frac{x_i-x}{h},\frac{y_i-y}{h}\Big)\Big]\Big|
\le 2\|K\|_{\infty},
\end{equation}
in the light of (\ref{K:Lq}).

We combine Fubini-Tonelli Theorem, (\ref{K:Lq}), (\ref{eYi22}) with (\ref{esup}) and since $f$ is a pdf, we deduce that
\begin{align}\label{S1}
S_1&\leq 2^{p-2}n^{1-p}h^{-p(d_1+d_2)}\|K\|_{\infty}^{p-2}\iint_{\bR^{d_1}\times\bR^{d_2}} \bE_f\big(\eta_i^2(x,y)\big)dxdy\\ \nonumber
&\leq 2^{p-2}n^{1-p}h^{-p(d_1+d_2)}\|K\|_{\infty}^{p-2}\iint_{\bR^{d_1}\times\bR^{d_2}}\iint_{\bR^{d_1}\times\bR^{d_2}} K^2\Big(\frac{z-x}{h},\frac{w-y}{h}\Big)f(z,w) dzdwdxdy\\ \nonumber
&=2^{p-2}n^{1-p}h^{(1-p)(d_1+d_2)}\|K\|_{\infty}^{p-2}\iint_{\bR^{d_1}\times\bR^{d_2}} K^2(x,y)dxdy\\ \nonumber
&\leq 2^{p-2}\|K\|_{\infty}^{p-2}\|K\|_2^2\big(nh^{d_1+d_2}\big)^{-p/2},
\end{align}
where fthe last inequality is due to the fact that $p\ge2$ and $nh^{d_1+d_2}\ge1$.

{\vspace{0.3cm}}

\textit{Estimation of $S_2$.} 

Inequality (\ref{eYi22}) implies
\begin{align}\nonumber
S_2&\leq n^{-p/2}h^{-p(d_1+d_2)}\iint_{\bR^{d_1}\times\bR^{d_2}} \Big(\mathbb{E}_f\big(\eta_i^2(x,y)\big)\Big)^{p/2}dxdy\\ \nonumber
&\leq  n^{-p/2}h^{-p(d_1+d_2)}\iint_{\bR^{d_1}\times\bR^{d_2}}\Big(\iint_{\bR^{d_1}\times\bR^{d_2}} K^2\Big(\frac{z-x}{h},\frac{w-y}{h}\Big)f(z,w) dzdw\Big)^{p/2}dxdy\\ \nonumber
&=:n^{-p/2}h^{-p(d_1+d_2)}\|Tf\|^{p/2}_{p/2},
\end{align}
where $T$ is the integral operator with kernel
\begin{equation}\label{T0}
T((z,w),(x,y)):=K^2\Big(\frac{z-x}{h},\frac{w-y}{h}\Big),\;\quad (z,w), (x,y)\in\bR^{d_1}\times\bR^{d_2}.
\end{equation}
A change of variable shows that 
\begin{equation}\label{T1}
\|T((z,w),\cdot)\|_1=\|T(\cdot,(x,y))\|_1=h^{d_1+d_2}\|K\|_2^2.
\end{equation}
Therefore, \cite[Theorem 6.36]{Folland} shows that
$$\|Tf\|_{p/2}^{p/2} \leq h^{(d_1+d_2)p/2}\|K\|_2^p \|f\|_{p/2}^{p/2}.
$$
In conclusion, we have
\begin{align}\label{S2}
S_2&\le \|K\|_2^p\big(nh^{d_1+d_2}\big)^{-p/2}\|f\|_{p/2}^{p/2}
\le \|K\|_2^p\big(nh^{d_1+d_2}\big)^{-p/2}\|f\|_{p}^{\frac{p(p-2)}{2(p-1)}}
\\
\label{S222}
&\le c\big(nh^{d_1+d_2}\big)^{-p/2}\|f\|_{W^{(s_1,s_2)}_p}^{\frac{p(p-2)}{2(p-1)}},
\end{align}
where we used Riesz-Thorin's Theorem, since in this case $1<p/2<p$.

Finally, by combining (\ref{Tr}), (\ref{Biasupperbound}), (\ref{S}), (\ref{S1}) and (\ref{S2}) we arrive at (\ref{upper1}) and the proof of Theorem \ref{TH:UBpbig} is completed.

\subsubsection{Proof of Theorem \ref{TH:UBpsmall2}} Let $1\le p<2$. We follow (\ref{S}), we use the convexity inequality (\ref{Rosenthal2}) and (\ref{eYi22}) to derive
\begin{align}\label{I}
S&=\mathbb{E}_f\big\|\hat{f}_n-\bE_f[\hat{f}_n]\big\|_p^p=\mathbb{E}_f\Big\|\frac{1}{n}\frac{1}{h^{d_1+d_2}}\sum\limits_{i=1}^{n}\eta_i(x,y)\Big\|_p^p
\\ \nonumber
&=n^{-p}h^{-p(d_1+d_2)}\iint_{\bR^{d_1}\times\bR^{d_2}}\mathbb{E}_f\Big(\Big|\sum_{i=1}^n \eta_i(x,y)\Big|^p\Big) dxdy
\\ \nonumber
&\leq
n^{-p}h^{-p(d_1+d_2)}\iint_{\bR^{d_1}\times\bR^{d_2}}\Big(\sum_{i=1}^n\mathbb{E}_f\big( \eta_i^2(x,y)\big)\Big)^{p/2}dxdy\\ \nonumber
&\leq  n^{-p/2}h^{-p(d_1+d_2)}\iint_{\bR^{d_1}\times\bR^{d_2}}\Big(\iint_{B((x_0,y_0),R)}K^2\Big(\frac{z-x}{h},\frac{w-y}{h}\Big)f(z,w) dzdw\Big)^{p/2}dxdy\\ \nonumber
&=:n^{-p/2}h^{-p(d_1+d_2)}I,
\end{align}
since $\supp f\subset B((x_0,y_0),R)$.

{\vspace{0.3cm}}

We separate the integral $I$ in the sum of the integrals over the ball $B:=B((x_0,y_0),2R^{*})$ and its complement $B^c:=\bR^{d_1}\times\bR^{d_2}\setminus B$, where $R^{*}:=\max(R,1/2)$. We denote by
$$I_1:=\iint_{B_2}\Bigg (\iint_{B((x_0,y_0),R)}K^2\Big(\frac{z-x}{h},\frac{w-y}{h}\Big)f(z,w) dzdw\Bigg)^{p/2}dxdy
$$

$$I_2:=\iint_{B_2^c}\Bigg(\iint_{B((x_0,y_0),R)}K^2\Big(\frac{z-x}{h},\frac{w-y}{h}\Bigg )f(z,w) dzdw\Big)^{p/2}dxdy$$
and then
\begin{equation}
\label{thbub22}
I=I_1+I_2.
\end{equation}

{\vspace{0.3cm}}

\textit{Estimation of }$I_1$. Since $p<2$, there exists a $1<q<\infty$ such that $\frac{p}{2}+\frac{1}{q}=1$. By H\"{o}lder's inequality and Fubini's theorem we obtain
\begin{align}
I_1&\le|B_2|^{\frac{1}{q}}\Big(\iint_{\bR^{d_1}\times\bR^{d_2}}\iint_{\bR^{d_1}\times\bR^{d_2}}K^2\Big(\frac{z-x}{h},\frac{w-y}{h}\Big)f(z,w) dzdwdxdy\Big)^{p/2}
\nonumber
\\
&=|B_2|^{\frac{1}{q}}\Big(\iint_{\bR^{d_1}\times\bR^{d_2}}f(z,w) dzdw\iint_{\bR^{d_1}\times\bR^{d_2}}K^2\big(x',y'\big)h^{(d_1+d_2)}dx'dy'\Big)^{p/2}
\nonumber
\\
&=v(d_1+d_2)^{\frac{1}{q}}\big(2R^{*}\big)^{\frac{d_1+d_2}{q}}\|K\|_2^{p}h^{(d_1+d_2)p/2},
\label{thub23}
\end{align}
where we applied a chang of variables, and we used the fact that $f$ is a pdf. Moreover, we call $v(d):=\frac{\pi^{d/2}}{\Gamma(1+d/2)}$ the volume of the unit ball on $\bR^d$.

{\vspace{0.3cm}}

\textit{Estimation of }$I_2$. Since $K$ is assumed to be compactly supported, there exists a positive number $\alpha>0$ such that
$$\supp K\subset\{(x,y)\in\bR^{d_1}\times\bR^{d_2}:\;|(x,y)|<\alpha\}=:B_{\alpha}.$$

Let $(z,w)\in B((x_0,y_0),R)$, $(x,y)\in B_2^{c}$ and $\big(\frac{z-x}{h},\frac{w-y}{h}\big)\in\supp K\subset B_{\alpha}$. Recall that $R^{*}\ge1/2$. By triangle inequality, we derive
\begin{align*}
1+|(x,y)-(x_0,y_0)|&\le 2R^{*}+|(x,y)-(x_0,y_0)|<2|(x,y)-(x_0,y_0)|
\\
&\le2(|(x,y)-(z,w)|+|(z,w)-(x_0,y_0)|)<2(|(x,y)-(z,w)|+R)
\\
&<4|(x-z,y-w)|<4\alpha h.
\end{align*}
Therefore,
\begin{equation*}
K^2\Big(\frac{z-x}{h},\frac{w-y}{h}\Big)\le\frac{\|K\|_{\infty}^2(4\alpha h)^{2(d_1+d_2+1)/p}}{(1+|(x,y)-(x_0,y_0)|)^{2(d_1+d_2+1)/p}},
\end{equation*}
for the above $(z,w)$ and $(x,y)$. Thus,

\begin{align}
I_2&\le\|K\|_{\infty}^p(4\alpha)^{d_1+d_2+1}h^{d_1+d_2+1}\Bigg (\iint_{\bR^{d_1}\times\bR^{d_2}}f(z,w) dzdw\Bigg )^{p/2}
\nonumber
\\
&\times\iint_{\bR^{d_1}\times\bR^{d_2}}\frac{dxdy}{(1+|(x,y)-(x_0,y_0)|)^{d_1+d_2+1}}
\nonumber
\\
&\le\tilde{v}(d_1+d_2)(4\alpha)^{d_1+d_2+1}\|K\|_{\infty}^p h^{(d_1+d_2)p/2},
\label{thub24}
\end{align}
since $f$ is a pdf, $h\le1$ and $p\le2$ and where we denoted $\tilde{v}(d):=\int_{\bR^d}(1+|x|)^{-d-1}dx$, for $d\in\bN$.

The proof of Theorem \ref{TH:UBpsmall2} is completed after combining (\ref{Tr}) and (\ref{Biasupperbound}) with (\ref{I})-(\ref{thub24}).

\subsubsection{Proof of Theorem \ref{TH:UBpsmall1}} Let $1\le p<2$. Since $K$ is compactly supported, there exists a constant $\alpha=\alpha_K>0$ such that $\supp K(x,y)\subset B_\alpha:=\{(x,y)\in\bR^{d_1}\times\bR^{d_2}:\;|(x,y)|<\alpha\}$. Let also $\omega\in L^{p/2}\cap L^{\infty}$ as in Definition \ref{D:classomega}. With no loss of generality, we assume $(x_0,y_0)=0_{(d_1,d_2)}$.

Following (\ref{I}) we have
\begin{equation}
\label{thub212}
S\le n^{-p/2}h^{-(d_1+d_2)p}I,
\end{equation}
where 
$$I=\iint_{\bR^{d_1}\times\bR^{d_2}}\Bigg (\iint_{B((x,y),\alpha h)}K^2\Big(\frac{z-x}{h},\frac{w-y}{h}\Bigg )f(z,w) dzdw\Big)^{p/2}dxdy
$$
where $B((x,y),\alpha h):=\{(z,w)\in\bR^{d_1}\times\bR^{d_2}:\;|(z,w)-(x,y)|<\alpha h\}$.

{\vspace{0.3cm}}

\textit{Estimation of} $I$. We express the domain of integration $\bR^{d_1}\times\bR^{d_2}$ as the union $\bR^{d_1}\times\bR^{d_2}=\bigcup_{k=0}^{\infty}A_k$, where
$A_0:=\{(x,y):|(x,y)|<\alpha h\}$ and $$A_k:=\{(x,y):k\alpha h\le|(x,y)|<(k+1)\alpha h\},\;\text{for every}\;k\ge1.$$
We set
\begin{equation}\label{Ik}
I_k:=\iint_{A_k}\Bigg (\iint_{B((x,y),rh)}K^2\Big(\frac{z-x}{h},\frac{w-y}{h}\Big)f(z,w) dzdw\Bigg )^{p/2}dxdy,
\end{equation}
for every $k\ge0$ and therefore
\begin{equation}\label{Iv}
I=\sum\limits_{k=0}^{\infty}I_{k}.
\end{equation}

{\vspace{0.3cm}}

\textit{Estimation of} $I_0$. Since $f(z,w)\le\omega(z,w)\le\|\omega\|_{\infty}<\infty$ we get
\begin{align}
\label{I0}
I_0&\le\|\omega\|_{\infty}^{p/2}h^{(d_1+d_2)p/2}\iint_{A_0}\Bigg (\iint_{\bR^{d_1}\times\bR^{d_2}}K^2(z,w) dzdw\Bigg )^{p/2}dxdy
\\
&=v(d_1+d_2)\alpha^{d_1+d_2}\|\omega\|_{\infty}^{p/2}\|K\|_2^{p}h^{(d_1+d_2)(1+p/2)}.
\nonumber
\end{align}

{\vspace{0.3cm}}

\textit{Estimation of} $I_k$, $k\ge1$. We denote by $\theta(\rho)$ the value of $\omega(z,w)$, for every $(z,w)\in\bR^{d_1}\times\bR^{d_2}$, with $|(z,w)|=\rho,$ $\rho\ge0$.

Let $k\ge1$, $(x,y)\in A_k$ and $(z,w)\in B((x,y),r\alpha h)$. By triangle inequality, we have $|(z,w)|\ge|(x,y)|-|(z,w)-(x,y)|\ge(k-1)\alpha h$. Therefore,
\begin{align}
\label{thetaomega}
f(z,w)&\le\omega(z,w)=\theta(|(z,w)|)\le\theta((k-1)\alpha h)
\\
&=\min\{\omega(x,y):(x,y)\in A_{k-1}\}.
\nonumber
\end{align}

By (\ref{thetaomega}) and after a change  of variable in (\ref{Ik}), we have
\begin{align}\label{Ik2}
I_k&\le h^{d_1+d_2}\|K\|_2^p\int_{A_k}\theta((k-1)\alpha h)^{p/2}dxdy
\\
&= h^{d_1+d_2}\|K\|_2^p\frac{|A_k|}{|A_{k-1}|}\int_{A_{k-1}}\omega(x,y)^{p/2}dxdy
\nonumber
\\
&\le 2^{d_1+d_2}h^{d_1+d_2}\|K\|_2^p\int_{A_{k-1}}\omega(x,y)^{p/2}dxdy.
\nonumber
\end{align}

{\vspace{0.3cm}}

\textit{Estimation of} $I$. By (\ref{Iv}), (\ref{I0}) and (\ref{Ik2}) we derive, after summation,
\begin{align}\label{I3}
I&=I_0+\sum_{k=1}^{\infty}I_k
\\
&\le v(d_1+d_2)\alpha^{d_1+d_2}\|\omega\|_{\infty}^{p/2}\|K\|_2^{p}h^{(d_1+d_2)(1+p/2)}
\nonumber\\
&+2^{d_1+d_2}\|K\|_2^{p}h^{(d_1+d_2)p/2}\sum_{k=1}^{\infty}\int_{A_{k-1}}\omega(x,y)^{p/2}dxdy
\nonumber
\\
&= v(d_1+d_2)\alpha^{d_1+d_2}\|\omega\|_{\infty}^{p/2}\|K\|_2^{p}h^{(d_1+d_2)(1+p/2)}
\nonumber
\\
&+2^{d_1+d_2}\|K\|_2^{p}h^{(d_1+d_2)p/2}\|\omega\|_{p/2}^{p/2}.
\nonumber
\end{align}

We now combine (\ref{I}) with (\ref{thub212}) and (\ref{I3}) to conclude
\begin{equation}
\label{Ssmall1}
S\le cn^{-p/2}h^{-(d_1+d_2)p/2},
\end{equation}
where $c=v(d_1+d_2)\alpha^{d_1+d_2}\|\omega\|_{\infty}^{p/2}\|K\|_2^{p}+2^{d_1+d_2}\|K\|_2^{p}\|\omega\|_{p/2}^{p/2}$, a positive constant depending on $d_1,d_2,K,p$ and $\omega$.

The combination of (\ref{Tr}), (\ref{Biasupperbound}) and (\ref{Ssmall1}) completes the proof.
 
\section{Proof of the lower bounds}\label{Prooflower} We now proceed to present the proof of Theorems \ref{th:lowerbound} and \ref{th:lowerboundnoncompact} which include the lower bounds. We shall need some preliminaries first.

\subsection{Auxiliary Results} We first state two crucial results needed for our approach. The first result is a variation of \cite[Theorem 2.4]{Tsybakov}, and can be found  in \cite{GL}.

\begin{lemma}\label{Laux1} Let $\mathbb{F}$ be a space of probability densities and assume that for any $n\in\bN$ sufficiently large, there exists a positive real number $\rho_n$ and a finite set $\{f_0\}\cup\{f_{\omega}:\;\omega\in\Omega_n\}\subset\mathbb{F}$ such that
\begin{equation}
\label{L11}
\big\|f_{\omega}-f_{\omega'}\big\|_{p}\ge 2\rho_n,\quad\text{for every}\;\;\omega,\omega'\in\Omega_n\cup\{0\},\;\text{with}\;\omega\neq\omega'.
\end{equation} 
\begin{equation}
\label{L12}
\limsup\limits_{n\rightarrow\infty}\frac{1}{|\Omega_n|^2}\sum_{\omega\in\Omega_n}\mathbb{E}_{f_0}\Big\{\frac{d\mathbb{P}_{f_{\omega}}}{d\mathbb{P}_{f_0}}(X^{(n)})\Big\}^2=:c_0<\infty.
\end{equation} 
Then for any $p\ge1$
\begin{equation}\label{lowerlemma1}
\liminf\limits_{n\rightarrow\infty}\inf_{\tilde{f}}\sup_{f\in \Omega_n}\rho_n^{-p}\mathbb{E}_f\|\tilde{f}-f\|_p^p\ge \big(\sqrt{c_0}+\sqrt{c_0+1}\big)^{-2},
\end{equation}
where the infimum is taken over all possible estimators.
\end{lemma} 
The celebrated Varshamov-Gilbert theorem is also needed, and is reported here for completeness of exposition. 
\begin{lemma}\label{Laux2}\cite{Tsybakov} Let $m\in\bN$ such that $m\ge8$ and let $\varrho_m$ be the Hamming distance in $\{0,1\}^m$, defined as
\begin{equation}
\varrho_m(a,b)=\sum_{j=1}^{m}\ONE\{a_j\neq b_j\},\quad\text{where}\;a=(a_1,\dots,a_m),b=(b_1,\dots,b_m)\in\{0,1\}^m.
\end{equation}
There exists a subset $P_m$ of $\{0,1\}^m$ such that $|P_m|\ge2^{m/8}$ and $\varrho_m(a,a')\ge\frac{m}{8}$ for every $a,a'\in P_m$ with $a\ne a'$.
\end{lemma} 

\subsection{Proof of Theorem \ref{th:lowerbound}}
We are going to prove that
\begin{equation}\label{lowerboundnew}
\liminf\limits_{n\rightarrow\infty}\Big(\frac{r_{*}^{D/S}}{n}\Big)^{-\frac{Sp}{2S+D}}\inf_{\tilde{f}}\sup_{f\in W^{(s_1,s_2)}_p (r,R)}\mathbb{E}_f\|\tilde{f}-f\|_p^p\ge c_0,
\end{equation}
for every $1\leq p< \infty$. Obviously this covers claim (ii) and claim (i) holds true since
\begin{eqnarray*}
c_0&\leq \liminf\limits_{n\rightarrow\infty}\Big(\frac{r_{*}^{D/S}}{n}\Big)^{-\frac{Sp}{2S+D}}\inf\limits_{\tilde{f}}\sup\limits_{f\in W^{(s_1,s_2)}_p (r,R)}\mathbb{E}_f\|\tilde{f}-f\|_p^p\\ 
&\leq \liminf\limits_{n\rightarrow\infty}\Big(\frac{r_{*}^{D/S}}{n}\Big)^{-\frac{Sp}{2S+D}}\inf\limits_{\tilde{f}}\sup\limits_{f\in W^{(s_1,s_2)}_p (r)}\mathbb{E}_f\|\tilde{f}-f\|_p^p.
\end{eqnarray*}

Our purpose is to apply Lemma \ref{Laux1} for $\mathbb{F}=W^{(s_1,s_2)}_{p}(r,R)$, for some $R>0$ sufficiently large and $\rho_n\sim n^{-S/(2S+D)}$; recall that $S=s_1+s_2$ and $D=d_1+d_2$.

We shall construct a family of functions $\{f_0,f_{\omega}:\;\omega\in\Omega\}\subset W^{(s_1,s_2)}_{p}(r,R)$ where the set $\Omega=\Omega_n$ will be chosen properly so that relations (\ref{L11}) and (\ref{L12}) to be satisfied.

{\vspace{0.3cm}}

We consider a function $k:\bR\rightarrow[0,\infty)$ such that:

(i) $k\in W^{s_1+s_2-1}_p(\bR)$, (ii) $\supp k\subset(-1,1)$ and (iii) $k$ is even\footnote{$k(-u)=k(u)$, for every $u\in\bR$}.
For example we can use the function $k(u)=\exp(-1/(1-u^2))\ONE_{(-1,1)}(u)$.

{\vspace{0.3cm}}

We define $\Lambda(u):=k(u)/\|k\|_1$, for every $u\in\bR$. It follows immediately that $\Lambda$ is a pdf supported in $(-1,1)$.

{\vspace{0.3cm}}

\textit{Step 1. Create} $f_0$: Consider a parameter $N>8$ and let the function
$$\bar{f}_0(x):=\prod_{\ell=1}^{D}\frac{1}{N}\big(\Lambda\ast\ONE_{[-\frac{N}{2},\frac{N}{2}]}\big)(x_{\ell}),\;\;\;\text{for every}\;x\in\bR^{D},$$
where $\ast$ denotes the convolution of the two functions; see (\ref{Conv}).

It can be namely verified that $\bar{f}_0$ is supported in $\big[-\frac{N}{2}-1,\frac{N}{2}+1\big]^D$. Further,
\begin{equation}
\label{thlb3}
\bar{f}_0 (x)=N^{-D},\;\;\forall x\in\Big[-\frac{N}{2}+1,\frac{N}{2}-1\Big]^D.
\end{equation}

We are going to show that $\bar{f}_0\in W^{(s_1,s_2)}_{p}$. By the construction of $\bar{f}_0$ we can express it in the form
\begin{equation}
\label{thlb33}
\bar{f}_0 (x)=N^{-D}\prod_{\ell=1}^{D}\tilde{\Lambda}(x_{\ell}),
\end{equation}
where $\tilde{\Lambda}(u):=\bar{\Lambda}\big(u+\frac{N}{2}\big)-\bar{\Lambda}\big(u-\frac{N}{2}\big)$, $u\in\bR$ and $\bar{\Lambda}(u):=\int_{0}^{u}\Lambda(t)dt$, an antiderivative of $\Lambda$. Hence,
\begin{align}
\big\|\bar{f}_0\big\|_{W^{(s_1,s_2)}_{p}}&=\sum_{|\alpha_1|\le s_1, |\alpha_2|\le s_2}\big\|\partial_2^{\alpha_2}\partial_1^{\alpha_1}\bar{f}_0\big\|_{p}
\nonumber\\
&=N^{-D}\sum_{|\alpha_1|\le s_1, |\alpha_2|\le s_2}\big\|\partial_2^{\alpha_2}\partial_1^{\alpha_1}\prod_{\ell=1}^{D}\tilde{\Lambda}(x_{\ell})\big\|_{p}
\nonumber\\
&\le N^{-D}\big\|\tilde{\Lambda}\big\|_{W^{s_1+s_2}_{p}(\bR)}^{D}.
\label{thlb4}
\end{align}

Apparently, $\|\tilde{\Lambda}^{(s)}\big\|_{p}^{p}=2\|\Lambda^{(s-1)}\big\|_{p}^{p}$, for every $s\ge1$.  Moreover, by Young's inequality (\ref{Young}), $\|\tilde{\Lambda}\big\|_{p}\le N\|\Lambda\|_p$. Since $N>2$ and $p\ge1$, we conclude
\begin{align}
\big\|\tilde{\Lambda}\big\|_{W^{s_1+s_2}_{p}(\bR)}&=\sum_{s\le s_1+s_2}\big\|\tilde{\Lambda}^{(s)}\big\|_{p}
\le N\|\Lambda\|_p+2^{1/p}\|\Lambda\|_{W^{s_1+s_2-1}_{p}(\bR)}
\nonumber\\
&\le 2N\frac{\|k\|_{W^{s_1+s_2-1}_{p}(\bR)}}{\|k\|_1},
\label{thlb41}
\end{align}
by definition of $\Lambda$. Finally by (\ref{thlb4}) and (\ref{thlb41}) we have
\begin{equation}
\label{thlb5}
\big\|\bar{f}_0\big\|_{W^{(s_1,s_2)}_{p}}\le\Big(\frac{2\|k\|_{W^{s_1+s_2-1}_{p}(\bR)}}{\|k\|_1}\Big)^{D}=:C_0<\infty
\end{equation}
and hence $\bar{f}_0\in W^{(s_1,s_2)}_{p}(C_0)$.

{\vspace{0.3cm}}

We are now ready to define $$f_0(x):=\kappa^{D}\bar{f}_0 (\kappa x),\;\;\text{for}\;\;0<\kappa\le1.$$ Then the pdf $f_0$ is clearly supported over
\begin{equation}
\label{thlbf0support}
\supp f_0(x)\subset\Big[-\frac{N+2}{2\kappa},\frac{N+2}{2\kappa}\Big]^D.
\end{equation}

We also note that  (\ref{thlb3}) implies
\begin{equation}
\label{thlbf0constant}
f_0(x)=\Big(\frac{\kappa}{N}\Big)^D,\quad\text{for every}\;x\in\Big[\frac{-N+2}{2\kappa},\frac{N-2}{2\kappa}\Big]^D.
\end{equation}

We now provide an estimate for the mixed Sobolev norm of $f_0$. Let us distinguish the cases $p=1$ and $p>1$:

{\vspace{0.3cm}}

For $p=1$, we recall that $r>1$. 

We derive
$$\big\|f_0\big\|_{W^{(s_1,s_2)}_{1}}=\|f_0\|_1+\sum_{|\alpha_1|\le s_1, |\alpha_2|\le s_2, 1\le|\alpha_1|+|\alpha_2|}\kappa^{|\alpha_1|+|\alpha_2|}\big\|\partial_2^{\alpha_2}\partial_1^{\alpha_1}\bar{f}_0\big\|_{1}\le 1+\kappa C_0,$$
in the light of (\ref{thlb5}) and since $\kappa\le1$.

Let $0<\varepsilon<1$ be such that $\varepsilon r>1$. Then, $f_0\in W^{(s_1,s_2)}_{p}(\varepsilon r)$, for $\kappa$ sufficiently small; $\kappa\le\frac{\varepsilon r-1}{C_0}$.

{\vspace{0.3cm}}

For $p>1$ it holds
$$\big\|f_0\big\|_{W^{(s_1,s_2)}_{p}}=\kappa^{D(1-\frac{1}{p})}\sum_{|\alpha_1|\le s_1, |\alpha_2|\le s_2}\kappa^{|\alpha_1|+|\alpha_2|}\big\|\partial_2^{\alpha_2}\partial_1^{\alpha_1}\bar{f}_0\big\|_{p}.$$
Then, $f_0\in W^{(s_1,s_2)}_{p}(\varepsilon r)$, for $\kappa$ be such that $\kappa\leq(\varepsilon rC_0^{-1})^{1/D(1-\frac{1}{p})}$.
Therefore, in every case we have, under the correct choice of $\kappa$,
\begin{equation}\label{th1fofinal}
f_0\in W^{(s_1,s_2)}_{p}(\varepsilon r, R),
\end{equation}
for some $0<\varepsilon<1$ and $R$ large enough, namely by (\ref{thlbf0support}) it is enough $R>\frac{\sqrt{D}}{\kappa} \big(\frac{N}{2}+1\big)$.

{\vspace{0.3cm}}

\textit{Step 2. Create the index set} $\Omega$: We start by defining a parameter $\sigma=\sigma(n)\rightarrow0$, as $n\rightarrow\infty$ with $\sigma<\min\big(1,\frac{1}{20\kappa}\big)$. In the sequel we introduce the parameter
\begin{equation}\label{thlbM}
 M:=\frac{N}{20\kappa\sigma}
\end{equation} 
and without loss of generality we assume that $M$ is an integer.

Let $\mathcal{M}=\{1,\dots,M\}^{D}$. By (\ref{thlbM}) is holds that $|\mathcal{M}|=M^{D}\ge8$. Hence, in application of Lemma \ref{Laux2}, we can find a set $\Omega\subset\{0,1\}^{|\mathcal{M}|}$ such that:
\begin{equation}
\label{thlb11}
|\Omega|\ge2^{|\mathcal{M}|/8}
\end{equation}
and for every $\omega,\omega'\in\Omega$ with $\omega\neq\omega'$
\begin{equation}
\label{thlb12}
\varrho_{|\mathcal{M}|}(\omega,\omega')\ge\frac{|\mathcal{M}|}{8},
\end{equation}
where $\varrho_{|\mathcal{M}|}$ is the Hamming distance in $\{0,1\}^{|\mathcal{M}|}$.

{\vspace{0.3cm}}

\textit{Step 3. Create} $f_{\omega}$: Using the index set $\Omega$ introduced in Step {\em 2.}, we will define the functions $f_{\omega}$. 

We first need to introduce a number of auxiliary functions.

Let us begin by defining the function $$g(t):=\Lambda\ast\big(\ONE_{[0,1]}-\ONE_{[-1,0]}\big)(t),\quad\text{for every}\;\;t\in\bR.$$ 
Then, by the properties of $\Lambda$, it turns out that:

(i) $g\in W^{s_1+s_2}_{p}(\bR)$, (ii) $\int_{\bR}g(y)dy=0$, (iii) $\supp g\subset[-2,2]$ and (iv) $\|g\|_{\infty}\le1$, thanks to Young's inequality (\ref{Young}).

Define now $$\xi_j=-\frac{N-4}{4\kappa}+8j\sigma,\quad\text{for}\;\;j=1,\dots,M.$$

For $m=(m_1,\dots,m_{D})\in\mathcal{M}$, set
$$G_m(x):=\prod_{j=1}^{D}g\Big(\frac{x_j-\xi_{m_j}}{\sigma}\Big),\quad\text{for}\;\;x\in\bR^D.$$

By the properties of $g$ we get
\begin{equation}
\label{thlb8}
\|G_m\|_{\infty}\le1,
\end{equation}
\begin{equation}
\label{thlb9}
\supp G_m\subset\Pi_m:=[\xi_{m_1}-3\sigma,\xi_{m_1}+3\sigma]\times\cdots[\xi_{m_D}-3\sigma,\xi_{m_D}+3\sigma].
\end{equation}
Note that
\begin{equation}
\label{thlb10}
\Pi_m\cap\Pi_{m'}\neq\emptyset,\quad \text{for every}\;m,m'\in\mathcal{M},\;\text{with}\;m\neq m'.
\end{equation}
Moreover
\begin{equation}
\label{thlb100}
\int_{\bR^D}G_m(x)dx=0,\quad \text{for every}\;m\in\mathcal{M}
\end{equation}
and by simple changes of variable,
\begin{equation}
\label{thlb1000}
\int_{\bR^D}|G_m(x)|^pdx=\sigma^{D}\|g\|_{p}^{pD},\quad \text{for every}\;m\in\mathcal{M}.
\end{equation}

Let the mapping $\pi:\mathcal{M}\rightarrow\{1,\dots,|\mathcal{M}|\}$ be such that
$$\pi(m):=\sum_{j=1}^{D}\big(m_{j}-1\big)M^{D-j}+m_{D}.$$
Note that $\pi$ defines an enumeration of the set $\mathcal{M}$ and it is a bijection.

We continue by defining the auxiliary family of functions $\{F_{\omega}:\omega\in\Omega\}$ as
$$F_{\omega}(x):=A\sum_{m\in\mathcal{M}}\omega_{\pi(m)}G_m(x),\quad\text{for }x\in\bR^D,$$
where $\omega=(\omega_1,\dots,\omega_{|\mathcal{M}|})\in\Omega$ and $A$ a parameter that will be specified in the sequel.

From relations (\ref{thlb8})-(\ref{thlb10}) we have
\begin{equation}
\label{thlb1200}\|F_{\omega}\|_{\infty}\le A\quad\text{for every}\;\omega\in\Omega.
\end{equation} 
By (\ref{thlb100}) it turns out that
\begin{equation}
\label{thlb120}
\int_{\bR^D}F_{\omega}(x)dx=0,\quad\text{for every}\;\omega\in\Omega
\end{equation}
and thanks to (\ref{thlb9})
\begin{equation}
\label{thlb16}
\supp F_{\omega}\subset\Big[-\frac{N-4}{4\kappa},\frac{N+4}{4\kappa}\Big]^{D}.
\end{equation}

We check now under which assumptions 
\begin{equation}
\label{thlbFomeganorm}
\|F_{\omega}\|_{W^{(s_1,s_2)}_{p}}\le r(1-\varepsilon).
\end{equation}

By the definition of $F_{\omega}$ and in the light of (\ref{thlb9}) and (\ref{thlb10}) we have for every $\alpha_1=(\alpha_{11},\dots,\alpha_{1d_1})\in\bN^{d_1}$ and $\alpha_2=(\alpha_{21},\dots,\alpha_{2d_2})\in\bN^{d_2}$, such that $|\alpha_1|\le s_1$ and $|\alpha_2|\le s_2$
\begin{align}
\big\|\partial^{\alpha_2}_2 \partial^{\alpha_1}_1 F_{\omega}\big\|_p^p&=\int_{\bR^D}\Big|\partial^{\alpha_2}_2\partial^{\alpha_1}_1\Big(A\sum_{m\in\mathcal{M}}\omega_{\pi(m)}G_m(x)\Big)\Big|^pdx
\nonumber
\\
&=A^p\sum_{m\in\mathcal{M}}\omega_{\pi(m)}\big\|\partial^{\alpha_2}_2\partial^{\alpha_1}_1 G_m\big\|_p^p.
\label{thlb13}
\end{align}

It is easy to verify that
\begin{align}
\big\|\partial^{\alpha_2}_2\partial^{\alpha_1}_1 G_m\big\|_p^p&=\prod_{j_1=1}^{d_1}\Big\|g\Big(\frac{x_{j_1}-\xi_{m_{j_1}}}{\sigma}\Big)^{(\alpha_{1j_1})}\Big\|_p^p
\prod_{j_2=1}^{d_2}\Big\|g\Big(\frac{x_{j_2}-\xi_{m_{j_2}}}{\sigma}\Big)^{(\alpha_{2j_2})}\Big\|_p^p
\nonumber
\\
&=\sigma^{D-p(|\alpha_1|+|\alpha_2|)}\prod_{j_1=1}^{d_1}\big\|g^{(\alpha_{1j_1})}\big\|_p^p\prod_{j_2=1}^{d_2}\big\|g^{(\alpha_{2j_2})}\big\|_p^p,
\label{thlb14}
\end{align}
where we used again a change of variable.

By combining (\ref{thlb13}) with (\ref{thlb14}) we conclude that
\begin{equation}
\label{thlbFw}
\|\partial^{\alpha_2}_2\partial^{\alpha_1}_1 F_{\omega}\big\|_p\le A|\mathcal{M}|^{1/p}\sigma^{\frac{D}{p}-(|\alpha_1|+|\alpha_2|)}\prod_{j_1=1}^{d_1}\big\|g^{(\alpha_{1j_1})}\big\|_p\prod_{j_2=1}^{d_2}\big\|g^{(\alpha_{2j_2})}\big\|_p.
\end{equation}

Given the above and since $\sigma<1$, we derive
\begin{align}
\|F_{\omega}\|_{W^{(s_1,s_2)}_p}&=\sum_{|\alpha_1|\le s_1, |\alpha_2|\le s_2}\|\partial_2^{\alpha_2}\partial_1^{\alpha_1}F_{\omega}\|_{p}
\nonumber
\\
&\le A|\mathcal{M}|^{1/p}\sigma^{\frac{D}{p}-S}\sum_{|\alpha_1|\le s_1, |\alpha_2|\le s_2}\prod_{j_1=1}^{d_1}\big\|g^{(\alpha_{1j_1})}\big\|_p\prod_{j_2=1}^{d_2}\big\|g^{(\alpha_{2j_2})}\big\|_p
\nonumber
\\
&=AM^{D/p}\sigma^{\frac{D}{p}-S}\|g\|_{W^{s_1+s_2}_{p}(\bR)}^{D},
\nonumber
\end{align}
where we recall that $S=s_1+s_2$. 

Then (\ref{thlbFomeganorm}) is guaranteed if we require
\begin{equation}
\label{thlb15}
A\Big(\frac{N}{20\kappa}\Big)^{\frac{D}{p}}\sigma^{-S}\|g\|_{W^{s_1+s_2}_{p}(\bR)}^{D}\le r(1-\varepsilon),
\end{equation}
thanks to (\ref{thlbM}).

We are now able to define, for any $\omega\in\Omega$,
$$f_{\omega}(x):=f_0 (x)+F_{\omega}(x),\quad\text{for every}\;x\in\bR^D.$$

Let us justify the properties of $f_{\omega}$'s.

{\vspace{0.3cm}}

(i) The functions $f_{\omega}$ are indeed pdf's:

Since $f_0$ is a pdf and thanks to (\ref{thlb120}) we obtain $\int_{\bR^D}f_{\omega}=1$.

We further need to ensure $f_{\omega}(x)\ge0$. By (\ref{thlb16}) and since $f_0$ is non-negative (as a pdf), it suffices to consider only the square $\big[-\frac{N-4}{4\kappa},\frac{N+4}{4\kappa}\big]^{D}$. But thanks to (\ref{thlbf0constant}) and (\ref{thlb1200}), it holds that $f_{\omega}(x)\ge0$ provided
\begin{equation}
\label{thlb18}
A\le \Big(\frac{\kappa}{N}\Big)^D.
\end{equation}

(ii) By (\ref{th1fofinal}) and (\ref{thlbFomeganorm}) we deduce
$$\|f_{\omega}\|_{W^{(s_1,s_2)}_{p}}\le \big\|f_0\big\|_{W^{(s_1,s_2)}_{p}}+\|F_{\omega}\|_{W^{(s_1,s_2)}_{p}}\le \varepsilon r+r(1-\varepsilon)=r.$$
and by (\ref{thlb16}) 
$$\supp f_{\omega} \subset \supp f_0,$$
so (\ref{th1fofinal}) gives
$$f_{\omega}\in W^{(s_1,s_2)}_{p}(r, R),$$
for $R$ large enough.

{\vspace{0.3cm}}

As a summary of all the above, we conclude that
$$\big\{f_0\big\}\cup\big\{f_{\omega}:\;\omega\in\Omega\big\}$$
is a finite set of pdf's contained in $W^{(s_1,s_2)}_{p}(r,R)$.

Thus, we are able to use Lemma \ref{Laux1} for $\mathbb{F}=W^{(s_1,s_2)}_{p}(r,R)$ and $\Omega$ being the set we introduced in Step {\em 2.} We must ensure under which assumptions, conditions (\ref{L11}) and (\ref{L12}) of Lemma \ref{Laux1} are fulfilled.

{\vspace{0.3cm}}

\textit{Step 4. Verifying condition (\ref{L11}):} By the definition of $f_{\omega}$'s and relations (\ref{thlb9}) and (\ref{thlb10}) we extract for every distinct $\omega,\omega'\in\Omega$
\begin{align*}
\big\|f_{\omega}-f_{\omega'}\big\|_p^p&=\big\|F_{\omega}-F_{\omega'}\big\|_p^p
\\
&=A^p\sum_{m\in\mathcal{M}}|\omega_{\pi(m)}-\omega'_{\pi(m)}|\int_{\Pi_m}|G_m(x)|^p dx
\\
&=A^p\varrho_{|\mathcal{M}|}(\omega,\omega')\sigma^{D}\|g\|_p^{pD}\ge A^p\frac{|\mathcal{M}|}{8}\sigma^{D}\|g\|_p^{pD}
\\
&=\frac{A^p}{8}\Big(\frac{N}{20\kappa}\Big)^D\|g\|_p^{pD},
\end{align*}
in the light of (\ref{thlb1000}), (\ref{thlb12}) and (\ref{thlbM}). Note that $f_0=f_{\omega_0}$, where $\omega_0=\{0,\dots,0\}$, so the above remains true for $\|f_{\omega}-f_{0}\|_p^p$ too.

Thus, we conclude that if $C_1:=\frac{1}{2}\|g\|_p^D\big(\frac{1}{20\kappa}\big)^{D/p}8^{-1/p}$, then condition (\ref{L11}) is fulfilled for
\begin{equation}
\label{thlbrn}
\rho_n:=C_1 AN^{D/p}.
\end{equation}

{\vspace{0.3cm}}

\textit{Step 5. Verifying condition (\ref{L12}):} It holds that
$$\frac{d\mathbb{P}_{f_{\omega}}}{d\mathbb{P}_{f_0}}\big(X^{(n)}\big)=\prod_{i=1}^{n}\frac{f_{\omega}(X_i)}{f_0 (X_i)}.$$
Moreover as $X_i$, $i=1,\dots,n$ are iid random variables for every $\omega\in\Omega$
\begin{align*}
\mathbb{E}_{f_0}\Big\{\prod_{i=1}^{n}\frac{f_{\omega}(X_i)}{f_0(X_i)}\Big\}^2&=\Big\{\int_{\bR^D}\frac{f_{\omega}^2(x)}{f_0 (x)}dx\Big\}^n
\\
&=\Big\{\int_{\bR^D}\frac{f_0^2 (x)+2f_0 (x)F_{\omega}(x)+F_{\omega}^2(x)}{f_0 (x)}dx\Big\}^n
\\
&=\Big\{1+\int_{\bR^D}\frac{F_{\omega}^2(x)}{f_0 (x)}dx\Big\}^n
\\
&=\big(1+\kappa^{-D}N^{D}\|F_{\omega}\|_2^2\big)^n,
\end{align*}
where for the third equality we used that $f_0$ is a pdf and (\ref{thlb120}), while for the last one, relations (\ref{thlb16}) and (\ref{thlbf0constant}) thanks to the the assumption $N>8$.

By (\ref{thlbFw}) and (\ref{thlbM}) it turns out that $\|F_{\omega}\|_2^2=A^2 \big(\frac{N}{20\kappa}\big)^{D}\|g\|_2^{2D}$ and hence for every $\omega\in\Omega$,
\begin{align*}
\mathbb{E}_{f_0}\Big\{\prod_{i=1}^{n}\frac{f_{\omega}(X_i)}{f_0(X_i)}\Big\}^2&=\big(1+\kappa^{-2D}N^{2D}20^{-D}A^2 \|g\|_2^{2D}\big)^n
\\
&\le\exp\big(n20^{-D}\kappa^{-2D}N^{2D}A^2\|g\|_2^{2D}\big),
\end{align*}
by the trivial inequality $1+x\le e^x$. Observe that the right hand side of the above is independent of $\omega$ thus,
\begin{align*}
\frac{1}{|\Omega|}\sum_{\omega\in\Omega}\mathbb{E}_{f_0}\Big\{\prod_{i=1}^{n}\frac{f_{\omega}(X_i)}{f_0(X_i)}\Big\}^2&\le\exp\big(C_2 nN^{2D}A^2\big),
\end{align*}
for $C_2:=20^{-D}\kappa^{-2D}\|g\|_2^{2D}$. Therefore if we require
\begin{equation*}
\ln|\Omega|\ge C_2 nN^{2D}A^2,
\end{equation*}
condition (\ref{L12}) is fulfilled for $c=1$. By (\ref{thlbM}) and (\ref{thlb11}) this is reduced to
\begin{equation}
\label{thlb19}
\Big(\frac{N}{20\kappa\sigma}\Big)^D \frac{\ln2}{8}\ge C_2 n N^{2D}A^2.
\end{equation}
{\vspace{0.3cm}}

\textit{Step 6. Choice of the parameters:} To summarize up, we constructed a family of functions $\{f_0,f_{\omega}:\;\omega\in\Omega\}\subset W^{(s_1,s_2)}_{p}(r,R)$ for which under the assumptions (\ref{thlb15}), (\ref{thlb18}) and (\ref{thlb19}) the assumptions of Lemma \ref{Laux1} are valid.

It only remains to choose the parameters $A$, $N$ and $\sigma$ so that (\ref{thlb15}), (\ref{thlb18}) and (\ref{thlb19}) to holds true and this ends the proof of the lower bound Theorem. 

We start by setting
\begin{equation*}
\varepsilon:=
\begin{cases}
\;\frac{r+1}{2r}\;,\;\;p=1\\
\frac{1}{2},\;\;p>1.
\end{cases}
\end{equation*}
Recall that the problem with the case $p=1$ is that we have to insure that $r>1$ and $\varepsilon r>1$, for being possible the corresponding Sobolev balls to contain densities. This fact causes the difference in the definition of $\varepsilon$.

Then the right hand side of (\ref{thlb15}) becomes
\begin{equation*}
r(1-\varepsilon)=\frac{r_*}{2}=
\begin{cases}
\;\frac{r-1}{2}\;,\;\;p=1\\
\;\frac{r}{2},\;\;p>1.
\end{cases}
\end{equation*}

Let also fix $N$ to be a constant with $N>8$. We set $C_3:=2\|g\|_{W^{S}_p(\bR)}^D\big(\frac{N}{20\kappa}\big)^{D/p}.$ We choose
\begin{equation}
\label{thlbsigma}
\sigma:=(C_3 A)^{\frac{1}{S}}r_*^{-\frac{1}{S}}=:C_4 A^{\frac{1}{S}}r_*^{-\frac{1}{S}}
\end{equation}
and then assumption (\ref{thlb15}) is fulfilled.

Given this choice of $\sigma$, for obtaining (\ref{thlb19}) it suffices to pick
\begin{equation}
\label{thlbA}
A:=C_5\Big(\frac{r_*^{D/S}}{n}\Big)^{\frac{S}{2S+D}},
\end{equation}
where
$$C_5:=\Big(\frac{\ln 2}{8 C_2 C_4^D N^D (20\kappa)^D}\Big)^{\frac{S}{2S+D}}.$$

By (\ref{thlbA}) replaced in (\ref{thlbrn}) and for $C_6:=C_1 C_5 N^{D/p}$ we arrive at
\begin{equation}
\label{thlbrn2}
\rho_n=C_6 \Big(\frac{r_* ^{D/S}}{n}\Big)^{\frac{S}{2S+D}}.
\end{equation}

By the choice of $A$ in (\ref{thlbA}), it turns out that $A\rightarrow0$, when $n\rightarrow\infty$. Therefore for sufficiently large values of $n$, condition (\ref{thlb18}) holds true. Moreover, by the choice of $\sigma$ in (\ref{thlbsigma}) we get that indeed $\sigma\rightarrow0$, when $n\rightarrow\infty$, consequently for sufficiently large values of $n$ it holds true that $\sigma<\min(1,1/(20\kappa))$, as assumed.

{\vspace{0.3cm}}

From all the above, Lemma 1 can be applied for
$\rho_n$ as in (\ref{thlbrn2}) and this completes the proof of Theorem \ref{th:lowerbound}.

\subsection{Proof of Theorem \ref{th:lowerboundnoncompact}} As we mentioned in Section \ref{ncpsmall} the lower bound has a different behaviour in the range $1\leq p<2$ if we do not assume that the density $f$ is compactly supported. 

The proof follows exactly the lines of the proof of Theorem \ref{th:lowerbound} with the only difference in the choice of the parameters.

The reason is that if we do not assume anymore that $N$ is a constant, we are able to achieve a greater lower bound for this case. Let us return to the Step {\em 6.} of Theorem \ref{th:lowerbound} and explain how the selection of our parameters $\sigma,\;A,\;N$ changes in this case.

Everything works mutatis mutandis, until the choice of $C_3$. We now fix $C_3':=2(20\kappa)^{-D/p}\|g\|_{W^S_p (\bR)}^D$ and choose
$$\sigma:=(C_3')^{1/S}A^{1/S}N^{Dp/S}r_{*}^{-1/S}=:C_4'A^{1/S}N^{Dp/S}r_{*}^{-1/S},$$
so that assumption (\ref{thlb15}) is fulfilled.

For this choice of $\sigma$, (\ref{thlb19}) is reduced to
\begin{equation}
\label{thlb20new}
A^{2+\frac{D}{S}}N^{D(1+\frac{Dp}{S})}r_{*}^{-\frac{D}{S}}\le n^{-1}C_5',
\end{equation}
for $C_5':=C_2^{-1}(C_4')^{-D}(20\kappa)^{-D}\frac{\ln 2}{8}$.

Let now $N^{D}=:C_6' A^{-1}$, such that $C_6'\le \kappa^D$ and hence relation (\ref{thlb18}) is fulfilled.

In order relation (\ref{thlb20new}) to be satisfied, it is enough to choose
\begin{equation}
\label{thlbAnew}
A:=C_7'n^{-\frac{pS}{pS+(p-1)D}}r_{*}^{\frac{pD}{pS+(p-1)D}},
\end{equation}
where $C_7':=\big(C'_5 C_6^{-{\frac{pS+D}{pS}}}\big)^{\frac{pS}{pS+(p-1)D}}$.

By replacing (\ref{thlbAnew}) in (\ref{thlbrn}) we obtain
\begin{equation}
\label{thlbrnnew}
\rho_n=C_8'\Big(\frac{r_{*}^{1/S}}{n}\Big)^{\frac{(p-1)S}{pS+(p-1)D}}.
\end{equation}

This time, $N\rightarrow\infty$ as $n\rightarrow\infty$. This is the choice that allows to elude the assumption of compactness of the support of the members of the family $\{f_0,f_{\omega}:\omega\in\Omega\}$. Note that, $\sigma\sim A^{\frac{1}{S}(1-\frac{1}{p})}$, hence for $p>1$, $\sigma\rightarrow 0$ as $n \rightarrow \infty$ because $A\rightarrow 0$ as $n \rightarrow \infty$. For $p=1$, we consider $C_6'$ to be small enough in order to guarantee that $\sigma\le\frac{1}{20\kappa}$, which is fundamental for the construction of the $f_{\omega}$'s, and the proof is complete.

\section{Final remarks} \label{Sec6} Let us close this paper with some discussions.

\subsection{Adaptivity}

One of the most interesting questions in the area of non-parametric estimation is to propose estimators that are \textit{adaptive}. This means that pdf's are provided estimates over a smoothness space $\bF^s_r$ while counting for the risk on the $L^p$-norm for a value of $p$ that might be different than $r$. For several results about adaptivity, the reader is referred to \cite{BKMP,CGKPP,DJKP,EY2,GL,GL3,HGPT,KePT,Rig}.

Studying the adaptivity issue for the setting proposed in this paper is a major challenge and will certainly be a fundamental target to be pursued by the authors in the future. 

Some necessary machinery for proposing adaptive estimators in the product setting we studied in this paper, it is not ready yet, but it will soon appear \cite{GKP}.

\subsection{Existence of kernels}\label{S:Kernels} The existence of kernels that belong to the class $\cK(s_1,s_2)$ is apparently not a trivial issue. The construction of meaningful examples can be based on tensor products of kernels that have been used for univariate functions in \cite{BH,Tsybakov}. 

In the Appendix of \cite{BH} or in Section 1 of \cite{Tsybakov} one can find examples of bounded kernels $\kappa:\bR\rightarrow\bR$ such that $\int_{\bR}\kappa=1$, $\int_{\bR}x^{\nu}\kappa(x)dx=0$, for $1\le\nu<s$, for $s\in\bN$ and $\int_{\bR}|x|^{s}|\kappa(x)|dx<\infty$. Such kernels use to be referred as kernels of order $s$.

Let $s_1,s_2\in\bN$ and $\kappa_1,\kappa_2$ be two kernels of order $s_1,s_2$ respectively. Then, it can be verified that the kernel $$K(x_1,\dots,x_{d_1},x_{d_1+1},\dots,x_{d_1+d_2}):=\kappa_1(x_1)\cdots\kappa_1(x_{d_1})\kappa_2(x_{d_1+1})\cdots\kappa_2(x_{d_1+d_2})$$
belongs to the class $\cK(s_1,s_2)$.

\subsection{Minimax estimation on classical Sobolev spaces \boldmath $W^s_p(\bR^d)$} 

Let $s\in\bN$. We denote by $\cK(s)$ the class of bounded Markov kernels $K:\bR^d\rightarrow\bR$ with vanishing moments $\int_{\bR^d}x^{\alpha}K(x)dx=0$, for every $1\le|\alpha|<s$ and such that $\int_{\bR^d}|x^{\alpha}||K(x)|dx<\infty$, for every $|\alpha|=s$.

Let us denote by $\cK_c(s)$ the subclass of $\cK(s)$ that contains compactly supported kernels.

The minimax problem on classical Sobolev spaces $W^{s}_p (\bR^d)$ takes the following form under the obvious modifications of the proofs of our theorems:
\begin{theorem}
\label{th:isotropic}
Let $r>0$, $1\le p<\infty$, $s\in\bN$, $K\in\cK(s)$ ($K\in\cK_c(s)$ for $p<2$) and $\hat{f}_n$ be the corresponding kernel density estimator defined in (\ref{KdeR}). Then
\begin{equation}\label{upperisotropicnew}
\sup_{f\in F^{s}_p}\mathbb{E}\|\hat{f}_n-f\|_p^p\le cn^{-sp/(2s+d)},
\end{equation}
where $F^{s}_p=W^s_p(r)$, when $p\ge2$ and $F^s_p=W^s_p(r,R)$, for some $R>0$, when $1\le p<2$. 
\end{theorem}
For the corresponding lower bounds, we have the following:
\begin{theorem}
\label{th:isotropic2}
Let $r_*$ as in relation (\ref{r*}), $1\le p<\infty$, $s\in\bN$, $F^s_p$ as in Theorem \ref{th:isotropic}. Then
\begin{equation}\label{lowerboundsmallpisotropic}
\liminf\limits_{n\rightarrow\infty}\Big(\frac{r_{*}^{d/s}}{n}\Big)^{-\frac{sp}{2s+d}}\inf_{\tilde{f}}\sup_{f\in F^s_p}\mathbb{E}\|\tilde{f}-f\|_p^p\ge c,
\end{equation}
where the infimum above is taken over all possible estimators $\tilde{f}$ and the constant $c$ is independent of $R$ and $r_*$.
\end{theorem}

As we mentioned earlier, the proofs of Theorems \ref{th:isotropic} and \ref{th:isotropic2} follow the lines of the results proved in 
Section \ref{Sec3} and are thus omitted.

Let us finally mention that theorems like those in Section \ref{ncpsmall} can be derived in a similar manner too.

\subsection{Minimax estimation on \boldmath $W^{(s_1,\dots,s_{\nu})}_p(\bR^{d_1}\times\cdots\times\bR^{d_{\nu}})$}\label{S:nproduct} The product spaces $W^{(s_1,s_2)}_p (\bR^{d_1}\times\bR^{d_2})$ that we studied, allow different level of smoothness on two sets of different variables. 

We prefer to present our contribution for this case just for simplicity. Our methods can be easily extended under very obvious modifications on the case of mixed regularity spaces like the following:
\begin{definition}\label{D:LSgen}
Let $\nu\ge2$, $s_1,\dots,s_{\nu}\in\mathbb{N}$, $1\le p<\infty$ and $f$ a function on $\bR^{d_1}\times\cdots\bR^{d_{\nu}}$. We say that $f$ belongs to  Sobolev space $W^{(s_1,\dots,s_{\nu})}_p(\bR^{d_1}\times\cdots\times\bR^{d_{\nu}})$, when
\begin{equation}\label{Sob3}\|f\|_{W^{(s_1,\dots,s_{\nu})}_p(\bR^{d_1}\times\cdots\times\bR^{d_{\nu}})}:=\sum_{|\alpha_i|\le s_i, i=1,\dots,\nu}\|\partial_{\nu}^{\alpha_{\nu}}\cdots\partial_1^{\alpha_1}f\|_{p}<\infty.
\end{equation}
\end{definition}
The results in Section \ref{Sec3} can be extended to this class of pdf's just by adapting the proofs we presented in Sections \ref{Sec4} and \ref{Prooflower}. Let us simply mention that this time the rate of convergence will be $n^{-pq_{\nu}}$, where
$$q_{\nu}:=\frac{s_1+\cdots s_{\nu}}{2(s_1+\cdots s_{\nu})+(d_1+\cdots d_{\nu})}.$$

\bibliographystyle{apalike}

\end{document}